\documentclass[final,onefignum,onetabnum]{siamart220329}

\usepackage{algorithm}
\usepackage[noend]{algpseudocode}
\usepackage{amsfonts}
\usepackage{amsmath}
\usepackage{amsopn}
\usepackage{amssymb}
\usepackage{cite}
\usepackage{enumitem}
\usepackage{float}
\usepackage{graphicx}
\usepackage[utf8]{inputenc}
\usepackage{latexsym}
\usepackage{lipsum}
\usepackage{mathtools}
\usepackage{tabularx}
\usepackage{xcolor}
\ifpdf
  \DeclareGraphicsExtensions{.eps,.pdf,.png,.jpg}
\else
  \DeclareGraphicsExtensions{.eps}
\fi

\newcommand{\R}{{\rm I\!R}}
\newcommand{\norm}[1]{\ensuremath{\left\|#1\right\|}}

\newsiamremark{remark}{Remark}
\newsiamremark{hypothesis}{Hypothesis}
\crefname{hypothesis}{Hypothesis}{Hypotheses}
\newsiamthm{claim}{Claim}
\newsiamremark{assumption}{Assumption}

\headers{EnKSGD}{B. Irwin and S. Reich}

\title{EnKSGD: A Class Of Preconditioned Black Box Optimization And Inversion Algorithms  \thanks{\funding{This work was supported by Deutsche Forschungsgemeinschaft (DFG) - Project-ID 318763901 - SFB1294.
The first author was also supported by the Natural Sciences and Engineering Research Council of Canada (NSERC) through a Postgraduate Scholarships - Doctoral (PGS-D) award. }}}

\author{Brian Irwin\thanks{Department of Earth, Ocean and Atmospheric Sciences, The University of British Columbia, Vancouver, British Columbia, Canada 
  (\email{birwin@eoas.ubc.ca}).}
\and Sebastian Reich\thanks{Department of Mathematics, University of Potsdam, Potsdam, Germany
  (\email{sebastian.reich@uni-potsdam.de}).}
}

\ifpdf
\hypersetup{
  pdftitle={EnKSGD: A Class Of Preconditioned Black Box Optimization And Inversion Algorithms},
  pdfauthor={B. Irwin and S. Reich}
}
\fi

\begin{document}

\maketitle

\begin{abstract}
In this paper, we introduce the Ensemble Kalman--Stein Gradient Descent (EnKSGD) class of algorithms. The EnKSGD class of algorithms builds on the ensemble Kalman filter (EnKF) line of work, applying techniques from sequential data assimilation to unconstrained optimization and parameter estimation problems. The essential idea is to exploit the EnKF as a black box (i.e. derivative-free, zeroth order) optimization tool if iterated to convergence. In this paper, we return to the foundations of the  EnKF as a sequential data assimilation technique, including its continuous-time and mean-field limits, with the goal of developing faster optimization algorithms suited to noisy black box optimization and inverse problems. The resulting EnKSGD class of algorithms can be designed to both maintain the desirable property of affine-invariance, and employ the well-known backtracking line search. Furthermore, EnKSGD algorithms are designed to not necessitate the subspace restriction property and variance collapse property of previous iterated EnKF approaches to optimization, as both these properties can be undesirable in an optimization context. EnKSGD also generalizes beyond the $L^{2}$ loss, and is thus applicable to a wider class of problems than the standard EnKF. Numerical experiments with both linear and nonlinear least squares problems, as well as maximum likelihood estimation, demonstrate the faster convergence of EnKSGD relative to alternative EnKF approaches to optimization.
\end{abstract}

\begin{keywords}
black box optimization, ensemble Kalman filter, affine-invariance, inverse problems, parameter estimation
\end{keywords}

\begin{MSCcodes}
65K10, 90C56, 65C35, 65C05, 62F10
\end{MSCcodes}

\section{Introduction}
\label{sec:intro}
This paper develops novel methods for the numerical solution of unconstrained optimization problems with the following structure
\begin{equation} \label{eq:general-inversion-objective}
\min_{\mathbf{x}} \bigg \{ \Phi \big ( \mathbf{x} \big ) \coloneqq \mathcal{D} \big ( \mathcal{G} ( \mathbf{x} ) \big ) + \alpha_{\mathbf{x}} \mathcal{R} ( \mathbf{x} ) + \alpha_{\mathbf{y}} \mathcal{T} \big ( \mathcal{G} ( \mathbf{x} ) \big ) \bigg \}
\end{equation}
where:
\smallskip
\begin{enumerate}
    \item 
    The objective function $\Phi(\mathbf{x}): \R^{N_{\mathbf{x}}} \mapsto \R$ is minimized with respect to the state variable $\mathbf{x} \in \R^{N_{\mathbf{x}}}$.
    \smallskip
    \item 
    The forward map $\mathcal{G}(\mathbf{x}): \R^{N_{\mathbf{x}}} \mapsto \R^{N_{\mathbf{y}}}$ transforms its input $\mathbf{x}$ from the state space $\R^{N_{\mathbf{x}}}$ into the observation space $\R^{N_{\mathbf{y}}}$. The forward map can be evaluated but its Jacobian matrix $D_{\mathbf{x}}\mathcal{G}(\mathbf{x}) \in \R^{N_{\mathbf{y}} \times N_{\mathbf{x}}}$ is unavailable.
    \smallskip    
    \item
    The loss function $\mathcal{D}(\mathbf{y}): \R^{N_{\mathbf{y}}} \mapsto \R$ transforms its input $\mathbf{y}$ from the observation space $\R^{N_{\mathbf{y}}}$ to a scalar real number. The loss function is convex and twice continuously differentiable, and its gradient $\nabla_{\mathbf{y}} \mathcal{D}(\mathbf{y}) \in \R^{N_{\mathbf{y}}}$ and Hessian $\nabla_{\mathbf{y}}^2 \mathcal{D}(\mathbf{y}) \in \R^{N_{\mathbf{y}} \times N_{\mathbf{y}}}$ are available.
    \smallskip
    \item
    The state space regularization function $\mathcal{R}(\mathbf{x}): \R^{N_{\mathbf{x}}} \mapsto \R$ transforms its input $\mathbf{x}$ from the state space $\R^{N_{\mathbf{x}}}$ to a scalar real number. The scalar real number $\alpha_{\mathbf{x}} \geq 0$ is the state space regularization parameter. Depending on the specific choice of state space regularization function, its gradient $\nabla_{\mathbf{x}} \mathcal{R}(\mathbf{x}) \in \R^{N_{\mathbf{x}}}$ and Hessian $\nabla_{\mathbf{x}}^2 \mathcal{R}(\mathbf{x}) \in \R^{N_{\mathbf{x}} \times N_{\mathbf{x}}}$ may or may not be available.
    \smallskip
    \item
    The observation space regularization function $\mathcal{T}(\mathbf{y}): \R^{N_{\mathbf{y}}} \mapsto \R$ transforms its input $\mathbf{y}$ from the observation space $\R^{N_{\mathbf{y}}}$ to a scalar real number. The scalar real number $\alpha_{\mathbf{y}} \geq 0$ is the observation space regularization parameter. Depending on the specific choice of observation space regularization function, $\nabla_{\mathbf{y}} \mathcal{T}(\mathbf{y}) \in \R^{N_{\mathbf{y}}}$ and $\nabla_{\mathbf{y}}^2 \mathcal{T}(\mathbf{y}) \in \R^{N_{\mathbf{y}} \times N_{\mathbf{y}}}$ may or may not be available. \\
\end{enumerate}
Optimization problems of the structure (\ref{eq:general-inversion-objective}) occur in diverse applications, such as:
\smallskip
\begin{enumerate}
    \item 
    Linear and nonlinear inverse problems, with and without common types of regularization, such as Tikhonov regularization.
    \smallskip
    \item
    Finding the roots of a function.
    \smallskip
    \item
    Empirical risk minimization (ERM) problems occurring in statistical estimation and machine learning. \\
\end{enumerate}
Our work is motivated by the fact that in practice, derivatives of the forward map $\mathcal{G}(\mathbf{x})$ are definitely not available, and derivatives of the regularization functions $\mathcal{R}(\mathbf{x})$ and $\mathcal{T}(\mathbf{y})$ may also not be available. This is because the user does not have access to a closed form expression for the function or because derivative computations may be too expensive to use in practice. As an example, the forward map $\mathcal{G}(\mathbf{x})$ may be a computer simulation tool used by engineers, and such simulation tools frequently do not provide a way for the user to obtain derivatives. As a result, in this case the forward map $\mathcal{G}(\mathbf{x})$ must be viewed as ``black box'' with unknown derivatives. Similarly, one may need to evaluate the output of a trained machine learning model, such as a neural network, to compute $\mathcal{R}(\mathbf{x})$ or $\mathcal{T}(\mathbf{y})$, but may not have access to the model's innner workings, and thus not have access to the derivatives of $\mathcal{R}(\mathbf{x})$ or $\mathcal{T}(\mathbf{y})$.

Furthermore, circumstances often arise where the forward map $\mathcal{G}(\mathbf{x})$ cannot be evaluated without stochastic noise $\eta(\mathbf{x})$ added. Mathematically, instead of $\mathcal{G}(\mathbf{x})$, one can evaluate $\mathcal{F}(\mathbf{x})$, the noisy forward map defined by $\mathcal{F}(\mathbf{x}) \coloneqq \mathcal{G}(\mathbf{x}) + \eta(\mathbf{x})$. A simple, yet relatively common, type of stochastic noise $\eta(\mathbf{x})$ is produced by drawing independent and identically distributed (IID) samples from a Gaussian distribution
\begin{equation} \label{eq:stochastic-noise-example}
    \eta(\mathbf{x}) \sim \mathcal{N}(\mathbf{0}, \sigma^2 \mathbf{I}_{N_{\mathbf{y}} \times N_{\mathbf{y}}}) \text{ , } \quad \sigma > 0 
\end{equation}
where $\mathbf{I}_{N_{\mathbf{y}} \times N_{\mathbf{y}}}$ denotes the $N_{\mathbf{y}} \times N_{\mathbf{y}}$ identity matrix. 

Finally, when dealing with multi-scale problems, it is desirable to employ methods which are invariant under affine changes of variables. That is, methods which are affine-invariant; a property shared by Newton's method but not by standard gradient descent methods. 

The characteristics mentioned above motivate us to develop methods based on ideas from the ensemble Kalman filter (EnKF) \cite{evensen-EnKF-original-paper-1994, evensen-EnKF-formulation-and-implementation-2003, evensen-EnKF-book-2009} and its variants. To provide the reader with important background material on the EnKF, in the following subsection we develop the essential ideas behind our proposed derivative-free optimization (DFO) method in the context of quadratic minimization problems. A central building block is provided by a continuous-time reformulation of the EnKF: the, so called, ensemble Kalman--Bucy equations \cite{bergemann-EnKF-localization-2010, reich-dynamical-2011}. This exposition is followed by a brief literature review examining the history and development of the EnKF, both in the context of Bayesian inference and optimization. This section closes with a summary of our key contributions and a layout of the remainder of the paper.

\subsection{The EnKF And Linear Inverse Problems}
\label{subsec:Kalman-Bucy-filter}
In this subsection, we explain the basic motivation for the proposed ensemble-based optimization algorithms for a linear inverse problem. Let us assume that the objective function $\Phi(\mathbf{x})$ in (\ref{eq:general-inversion-objective}) has regularization parameters $\alpha_{\mathbf{x}} = \alpha_{\mathbf{y}} = 0$, and takes the following quadratic form
\begin{equation} \label{eq:quadratic-objective-function-phi}
    \Phi(\mathbf{x}) = \frac{1}{2} \norm{\mathbf{G}\mathbf{x} - \mathbf{y}_{\rm obs}}_2^2
\end{equation}
for appropriate linear forward map $\mathcal{G}(\mathbf{x}) = \mathbf{G} \mathbf{x}$ with fixed matrix $\mathbf{G} \in \R^{N_{\mathbf{y}} \times N_{\mathbf{x}}}$, and quadratic loss function with fixed data vector $\mathbf{y}_{\rm obs} \in \R^{N_{\mathbf{y}}}$
\begin{equation} \label{eq:quadratic-loss-function}
    \mathcal{D}(\mathbf{y}) = \frac{1}{2} \norm{\mathbf{y} - \mathbf{y}_{\rm obs}}_2^2 \text{ . }
\end{equation}
We assume that $\Phi(\mathbf{x})$ is strongly convex and denote the unique minimizer by $\mathbf{x}^\dagger$. The continuous-time formulation of gradient descent is
\begin{equation} \label{eq:gradient-descent-continuous-time}
    \frac{{\rm d} \mathbf{x}(t)}{{\rm d}t} = -\nabla_\mathbf{x}\Phi \big ( \mathbf{x}(t) \big ) =
    -\mathbf{G}^{\rm T} \left( \mathbf{G} \mathbf{x}(t)-\mathbf{y}_{\rm obs}\right) \text{ . }
\end{equation}
It is well-known that gradient descent can converge slowly for linear multi-scale problems, as it lacks the property of affine-invariance. More precisely, an affine transformation with invertible matrix $\mathbf{A} \in \R^{N_{\mathbf{x}} \times N_{\mathbf{x}}}$ and vector $\mathbf{b} \in \R^{N_{\mathbf{x}}}$ of the form
\begin{equation} \label{eq:affine-transformation}
    \mathbf{x}(t) = \mathbf{A}\mathbf{\tilde x}(t) + \mathbf{b} 
\end{equation}
transforms (\ref{eq:gradient-descent-continuous-time}) into
\begin{equation} \label{eq:gradient-descent-continuous-time-affine-transformed}
    \frac{{\rm d} \mathbf{\tilde x}(t)}{{\rm d}t} = -\mathbf{A}^{-1}
    \mathbf{G}^{\rm T} \left( \mathbf{G}(\mathbf{A} \mathbf{\tilde x}(t) + \mathbf{b}) -\mathbf{y}_{\rm obs}\right) \text{ , }
\end{equation}
which differs from applying gradient descent to the transformed objective function
\begin{equation} \label{eq:quadratic-objective-function-phi-affine-transformed}
    \tilde{\Phi}(\mathbf{\tilde x}) = \frac{1}{2} \norm{\mathbf{G} (\mathbf{A} \mathbf{\tilde x} + \mathbf{b}) - \mathbf{y}_{\rm obs}}_2^2 \text{ . } 
\end{equation}
For that reason, one often uses Newton's method instead. The continuous-time formulation of Newton's method, which is an affine-invariant method, is
\begin{equation} \label{eq:Newton-continuous-time}
    \frac{{\rm d} \mathbf{x}(t)}{{\rm d}t}  =
    -\left(\mathbf{G}^{\rm T}
    \mathbf{G}\right)^{-1} \mathbf{G}^{\rm T} \left( \mathbf{G} \mathbf{x}(t)-\mathbf{y}_{\rm obs}\right) \text{ . }
\end{equation}
The affine-invariance of Newton's method is revealed by its equivalent formulation
\begin{equation} \label{eq:Newton-continuous-time-grad-Phi-formulation}
    \frac{\rm d}{{\rm d}t} \nabla_\mathbf{x} \Phi \big ( \mathbf{x}(t) \big )
    = - \nabla_\mathbf{x} \Phi \big ( \mathbf{x}(t) \big ) \text{ . }
\end{equation}
However, the affine-invariance property comes at a price; namely, the inversion of the symmetric positive-definite Hessian matrix $\mathbf{G}^{\rm T} \mathbf{G}$. 

We now develop an inversion-free approximation to Newton's method (\ref{eq:Newton-continuous-time}), which is based on continuous-time mean-field formulations of the EnKF \cite{bergemann-EnKF-localization-2010, reich-dynamical-2011}.  The first step is to treat the state variable $\mathbf{x}$ as a Gaussian random variable. We then consider the family of time-varying Gaussian distributions
\begin{equation} \label{eq:Phi-decaying-gaussian-family}
    \rho(\mathbf{x},t) \propto 
    \exp \big ( - t \Phi(\mathbf{x}) \big ) \rho_0(\mathbf{x})
\end{equation}
for $t \ge 0$ and given initial Gaussian distribution
\begin{equation} \label{eq:Phi-decaying-gaussian-initial}
    \rho_0 = \mathcal{N}({\mathbf{\bar x}}_0,\mathbf{P}_0) \text{ . }
\end{equation}
Letting $\delta(\mathbf{x})$ denote the Dirac delta function, it is easy to deduce that 
\begin{equation} \label{eq:Phi-decaying-gaussian-family-infinite-time}
    \lim_{t\to \infty} \rho(\mathbf{x},t) =\delta(\mathbf{x}-
    \mathbf{x}^\dagger) \text{ . }
\end{equation}
In other words, as $t \to \infty$, all the probability mass concentrates at the minimizer $\mathbf{x}^\dagger$. The associated, so called, Kalman--Bucy mean-field equation \cite{reich-dynamical-2011}
\begin{equation} \label{eq:Kalman-Bucy-mean-field}
    \frac{{\rm d} \mathbf{x}(t)}{{\rm d}t} =
    -\frac{1}{2} \mathbf{P}(t) \mathbf{G}^{\rm T}\left(
    \mathbf{G}\mathbf{x}(t)+ \mathbf{G}\mathbf{\bar x}(t)
    -2\mathbf{y}_{\rm obs}\right)
\end{equation}
satisfies
\begin{equation}
    \mathbf{x}(t) \sim \rho(\cdot,t) = \mathcal{N} \big (
    \mathbf{\bar x}(t),\mathbf{P}(t) \big ) \text{ . }
\end{equation}
Note that the Kalman--Bucy mean-field equation (\ref{eq:Kalman-Bucy-mean-field}) has a structure similar to Newton's method (\ref{eq:Newton-continuous-time}), with the key difference being that the inverse of the Hessian matrix $\mathbf{G}^{\rm T}\mathbf{G}$ in (\ref{eq:Newton-continuous-time}) is replaced by the covariance matrix $\mathbf{P}(t)$ in (\ref{eq:Kalman-Bucy-mean-field}). The Kalman--Bucy mean-field equation (\ref{eq:Kalman-Bucy-mean-field}) is affine-invariant \cite{pidstrigach-ensemble-transform-logistic-2022} but converges slowly to $\mathbf{x}^\dagger$ as
\begin{equation}
    \lim_{t \to \infty} \mathbf{P}(t) = \mathbf{0}
\end{equation}
at a rate proportional to $1/t$ \cite{schillings-EnKF-analysis-2017}.

Following the closely related work \cite{huang-derivative-free-Bayesian-inversion-2022} on Bayesian inference, we propose to instead consider the modified evolution equations with small scale parameter $\delta > 0$
\begin{equation} \label{eq:Kalman-Bucy-mean-field-modified}
    \frac{{\rm d} \mathbf{x}(t)}{{\rm d}t} = \frac{\delta}{2} \big ( \mathbf{x}(t)-\mathbf{\bar x}(t) \big )
    -\frac{1}{2} \mathbf{P}(t) \mathbf{G}^{\rm T}\left(
    \mathbf{G}\mathbf{x}(t)+ \mathbf{G}\mathbf{\bar x}(t)
    -2\mathbf{y}_{\rm obs} \right) \text{ . }
\end{equation}
The evolution in $\mathbf{x}(t)$ remains Gaussian with evolution equations for the mean
\begin{equation} \label{eq:Kalman-Bucy-mean-evolution}
    \frac{{\rm d} \mathbf{\bar x}(t)}{{\rm d}t} = 
    -\mathbf{P}(t) \mathbf{G}^{\rm T}\left(
    \mathbf{G}\mathbf{\bar{x}}(t)-\mathbf{y}_{\rm obs}\right) \text{ , }
\end{equation}
and the covariance matrix
\begin{equation} \label{eq:Kalman-Bucy-covariance-evolution}
    \frac{{\rm d} \mathbf{P}(t)}{{\rm d}t} = \delta \,\mathbf{P}(t)
    -\mathbf{P}(t) \mathbf{G}^{\rm T}
    \mathbf{G}\mathbf{P}(t) \text{ , }
\end{equation}
respectively. One finds that, asymptotically as $t \to \infty$,
\begin{equation}
    \mathbf{P}(t) \approx \mathbf{P}_\infty :=
    \delta \left(\mathbf{G}^{\rm T}
    \mathbf{G}\right)^{-1} \text{ . }
\end{equation}
Hence, the evolution in the mean asymptotically reduces to
\begin{equation}
    \frac{{\rm d} \mathbf{\bar x}(t)}{{\rm d}t} = 
    -\delta \left(\mathbf{G}^{\rm T}
    \mathbf{G}\right)^{-1} \mathbf{G}^{\rm T}\left(
    \mathbf{G}\mathbf{\bar x}(t)-\mathbf{y}_{\rm obs}\right) \text{ , }
\end{equation}
which is a time re-scaled version of Newton's method (\ref{eq:Newton-continuous-time}). Therefore, we finally consider the re-scaled mean-field equation
\begin{equation} \label{eq:two-scale-EnKBF}
    \frac{{\rm d} \mathbf{x}(t)}{{\rm d}t} = \frac{1}{2} \big ( \mathbf{x}(t)-\mathbf{\bar x}(t) \big )
    -\frac{1}{2\delta} \mathbf{P}(t) \mathbf{G}^{\rm T}\left(
    \mathbf{G}\mathbf{x}(t)+ \mathbf{G}\mathbf{\bar x}(t)
    -2\mathbf{y}_{\rm obs}\right) \text{ , }
\end{equation}
which is affine-invariant for any $\delta > 0$. Following our previous analysis, for the re-scaled mean-field equation (\ref{eq:two-scale-EnKBF}), one finds that $\mathbf{P}(t)$ converges rapidly to $\mathbf{P}_\infty$ for small values of $\delta$, while the evolution in the mean $\mathbf{\bar x}(t)$ mimics that of Newton's method for a quadratic objective function.

This paper studies generalization of the formulations (\ref{eq:Kalman-Bucy-mean-field-modified}) and (\ref{eq:two-scale-EnKBF}) to nonlinear minimization problems of the form (\ref{eq:general-inversion-objective}). In this context, observe that the gradient and Hessian of the quadratic choice of objective function $\Phi(\mathbf{x})$ given by (\ref{eq:quadratic-objective-function-phi}) are
\begin{equation}
    \nabla_{\mathbf{x}} \Phi \big ( \mathbf{\bar{x}}(t) \big ) = \mathbf{G}^{\rm T}  \left( \mathbf{G} \mathbf{\bar{x}}(t) - \mathbf{y}_{\rm obs} \right)
\end{equation}
and
\begin{equation}
    \nabla_{\mathbf{x}}^2 \Phi \big ( \mathbf{\bar x}(t) \big ) = \mathbf{G}^{\rm T} \mathbf{G}
\end{equation}
respectively. Thus, after rescaling by $1/\delta$, (\ref{eq:Kalman-Bucy-mean-evolution}) can be written as 
\begin{equation} \label{eq:Kalman-Bucy-mean-evolution-scaled}
    \frac{{\rm d} \mathbf{\bar{x}}(t)} {{\rm d} t} =  
    - \frac{1}{\delta} \mathbf{P}(t) \nabla_{\mathbf{x}} \Phi \big ( \mathbf{\bar{x}}(t) \big )
\end{equation}
and (\ref{eq:Kalman-Bucy-covariance-evolution}) as 
\begin{equation}  \label{eq:Kalman-Bucy-covariance-evolution-scaled}
    \frac{{\rm d} \mathbf{P}(t)}{{\rm d} t} =  \mathbf{P}(t) - \frac{1}{\delta}\mathbf{P}(t) \nabla_{\mathbf{x}}^2 \Phi \big ( \mathbf{\bar{x}}(t) \big ) \mathbf{P}(t) \text{ . }
\end{equation}
While (\ref{eq:Kalman-Bucy-covariance-evolution-scaled}) requires the Hessian of the objective function $\Phi(\mathbf{x})$, we propose a reformulation of Newton's method that not only avoids matrix inversions but also allows for derivative-free implementations while maintaining affine-invariance.

\subsection{Literature Review}
\label{subsec:literature-review}
A state-of-the-art introduction to unconstrained optimization methods can be found in \cite{nocedal-wright-numerical-optimization-2006}, including quasi-Newton methods, which avoid the explicit computation of Hessians and their inversion. Completely derivative-free methods have been recently surveyed in \cite{larson-dfo-review-2019}. In contrast to those derivative-free optimization methods, in this paper we build on extensions of the EnKF to optimization problems. Historically, the EnKF has been developed as a Monte Carlo based extension of the classical Kalman filter methodology \cite{kalman-new-linear-filtering-approach-1960} to nonlinear state and parameter estimation in the field of data assimilation \cite{evensen-EnKF-book-2009, law-data-assimilation-book-2015, reich-cotter-book-2015, asch-data-assimilation-book-2016}. An extensive survey of the EnKF in the context of Bayesian inference and optimization is provide by \cite{calvello-ensemble-Kalman-mean-field-perspective-ArXiv-2022}. Alternative extensions of the classical Kalman filter to nonlinear and non-Gaussian data assimilation problems include the extended Kalman filter (ExKF) \cite{bellantoni-Kalman-Schmidt-square-root-1967} and the unscented Kalman filter (UKF) \cite{julier-UKF-1997}.

\subsubsection{EnKF For Bayesian Inference}
\label{subsubsec:EnKF-literature-review}
Since its introduction in 1994 \cite{evensen-EnKF-original-paper-1994}, the EnKF has established itself as a favourite data assimilation algorithm of geoscience practitioners due to its wide applicability, ease of implementation, robustness, and computational efficiency. From a Bayesian perspective, the analysis step of the Kalman filter produces the posterior distribution of the state conditioned on observations. In the case of linear Gaussian dynamics, the EnKF converges to the true posterior distribution as the ensemble size approaches infinity (i.e. in the mean-field limit) \cite{legland-EnKF-asymptotics-2009, kwiatkowski-square-root-EnKF-convergence-2015}. However, in the case of nonlinear dynamics, the paper \cite{ernst-EnKF-analysis-2015} establishes that the EnKF does not converge to the true posterior distribution as the ensemble size approaches infinity. Similarly, \cite{law-deterministic-EnKF-2016} shows that the mean-field EnKF provides the optimal linear estimator of the conditional mean, but only produces the true posterior distribution in the case of linear Gaussian dynamics. Further analysis of the EnKF focusing on data assimilation in the large ensemble limit can be found in \cite{li-EnKF-numerical-properties-2008}.

In practice, the EnKF is most frequently used with ensembles of a fixed finite size much smaller than the state space dimension \cite{bergemann-EnKF-localization-2010, bergemann-mollified-EnKF-2010}. However, until the past decade, there was relatively little research analyzing the behaviour of the EnKF for fixed size ensembles. Recent research analyzing the behaviour of the EnKF for fixed size ensembles in the context of data assimilation includes \cite{kelly-EnKF-well-posedness-2014, tong-nonlinear-stability-EnKF-2016}. The use of finite size ensembles introduces errors that produce spurious correlations. Over time, the spurious correlations may lead to an undesirable reduction of the ensemble variance, with the ensemble ultimately collapsing to its mean. Several strategies have been developed to combat ensemble collapse, including variance inflation, localization, and sampling based strategies \cite{evensen-EnKF-book-2009}.   

More recently, the EnKF has also been combined with sampling methods for Bayesian inference. In \cite{garbuno-inigo-EKS-2020}, the authors introduce  the ensemble Kalman sampler (EKS). The EKS is a derivative-free approximate sampler for the posterior distribution. Unlike the EnKF, the EKS is designed to transform arbitrary samples at $t = 0$ into posterior samples as $t \rightarrow \infty$. EKS is further discussed in \cite{nusken-EKS-note-2019, ding-EKS-analysis-2021}. Furthermore, building again on Gaussian approximations, a very efficient Bayesian inference method has been proposed in \cite{huang-derivative-free-Bayesian-inversion-2022}. The essential idea is to combine the EnKF with an appropriate McKean--Vlasov type dynamical system which has a Gaussian approximation to the posterior as invariant distribution; see also \cite{pidstrigach-ensemble-transform-logistic-2022}. In this paper, we have used the same conceptional idea in deriving (\ref{eq:two-scale-EnKBF}) but now in the context of derivative-free optimization.

\subsubsection{EnKF For Optimization}
\label{subsubsec:EKI-literature-review}
Over the past decade, an alternative perspective viewing the EnKF as a derivative-free optimization algorithm has emerged from the work \cite{iglesias-EKI-2013}. This perspective has advanced a line of work \cite{iglesias-EKI-pde-constrained-2016, schillings-EnKF-analysis-2017, blomker-EKI-convergent-scheme-2018, albers-EKI-constrained-2019, ding-EKI-analysis-2021, blomker-EKI-convergence-analysis-2019, chada-EKI-parameterizations-2018, schillings-EKI-convergence-linear-noisy-2018, chada-EKI-tikhonov-2020} systematically developing the EnKF as a general tool for solving inverse problems. The essential idea is to keep iterating an EnKF implementation for fixed observations until the ensemble has collapsed onto the minimizer of the objective function. While most of the original work on the EnKF has focused on objective functions $\Phi(\mathbf{x})$ with a quadratic loss function (\ref{eq:quadratic-loss-function}), extensions to more general loss functions can be found in \cite{kovachki-EKI-machine-learning-2019, haber-never-look-back-EnKF-ArXiv-2018, pidstrigach-ensemble-transform-logistic-2022}. While the EnKF convergences at a rate of $1/t$ when used as a derivative-free optimization algorithm, combining the EnKF with appropriate dynamical updates leads to exponential convergence \cite{huang-iterated-Kalman-methodology-2022}.

We finally mention the class of consensus-based ensemble methods for optimization \cite{pinnau-consensus-optimization-2017}. These methods are also designed to perform derivative-free optimization but do not attempt to mimic a derivative-free version of Newton's method.

\subsection{Contributions}
\label{subsec:contributions}
Having introduced the background material in subsections~\ref{subsec:Kalman-Bucy-filter} and~\ref{subsec:literature-review}, we now outline the development of the main contribution of this paper: the Ensemble Kalman--Stein Gradient Descent (EnKSGD) class of algorithms. Section~\ref{sec:preconditioned-grad-flow} develops a particle-based approach to solving a preconditioned gradient flow for (\ref{eq:general-inversion-objective}). Section~\ref{sec:derivative-approximations} describes particle-based derivative approximations based in Stein's identity \cite{opper-bayesian-online-learning-1999, pidstrigach-ensemble-transform-logistic-2022}. Section~\ref{sec:EnKSGD-algorithms} combines the content of Sections~\ref{sec:preconditioned-grad-flow} and~\ref{sec:derivative-approximations} to produce the EnKSGD class of algorithms. Section~\ref{sec:numerical-experiments} presents the results of numerical experiments with minimization problems of the form (\ref{eq:general-inversion-objective}). Finally, Section~\ref{sec:conclusions} concludes the paper.

\section{Particle Based Approach To Solving Preconditioned Gradient Flow}
\label{sec:preconditioned-grad-flow}
Recall that we seek the solution of (\ref{eq:general-inversion-objective}). For a real, symmetric positive-definite matrix $\mathbf{B}^{-1} \big ( \mathbf{x}(t) \big ) \in \R^{N_{\mathbf{x}} \times N_{\mathbf{x}}}$, a preconditioned gradient flow for (\ref{eq:general-inversion-objective}) is given by
\begin{equation} \label{eq:general-inversion-preconditioned-gradient-flow}
    \frac{{\rm d} \mathbf{x}(t)} {{\rm d} t} = - \mathbf{B}^{-1} \big ( \mathbf{x}(t) \big ) \nabla_{\mathbf{x}} \Phi \big ( \mathbf{x}(t) \big ) \text{ . }
\end{equation}
We are chiefly interested in the cases where $\mathbf{B} \big ( \mathbf{x}(t) \big )$ is the Hessian $\nabla_{\mathbf{x}}^2 \Phi \big ( \mathbf{x}(t) \big )$,  as in Newton's method (\ref{eq:Newton-continuous-time}), or a modification of the Hessian.

\subsection{Generalized Ensemble Kalman--Bucy Mean-Field Equations}
\label{subsec:two-time-scale-approach}
Following the probabilistic interpretation laid out in subsection \ref{subsec:Kalman-Bucy-filter} for linear inverse problems, we now extend the ensemble Kalman--Bucy formulation (\ref{eq:two-scale-EnKBF}) to general nonlinear inverse problems with objective function (\ref{eq:general-inversion-objective}).  We consider (\ref{eq:Kalman-Bucy-mean-evolution-scaled}) and a modification of (\ref{eq:Kalman-Bucy-covariance-evolution-scaled}). As we want $\mathbf{P}(t)$ to behave like $\delta \mathbf{B}^{-1} \big ( \mathbf{\bar x}(t) \big )$, we modify (\ref{eq:Kalman-Bucy-covariance-evolution-scaled}) to obtain the modified Kalman--Bucy covariance equation given by
\begin{equation} \label{eq:modified-Kalman-Bucy-covariance-Phit}
    \frac{{\rm d} \mathbf{P}(t)}{{\rm d}t} = \mathbf{P}(t) - \frac{1}{\delta} \mathbf{P}(t) \mathbf{B} \big ( \mathbf{\bar x}(t) \big ) \mathbf{P}(t) \text{ . }
\end{equation}

We next analyse the implied dynamics in $\mathbf{P}(t)$ for constant $\mathbf{B} = \mathbf{B}(\mathbf{\bar x})$ with
$\mathbf{\bar x}$ held fixed. Using the following relationship 
\begin{equation} \label{eq:inverse-covariance-time-derivative}
    \frac{{\rm d} \mathbf{P}^{-1}(t)}{{\rm d} t} = - \mathbf{P}^{-1}(t) \frac{{\rm d} \mathbf{P}(t)}{{\rm d} t} \mathbf{P}^{-1}(t)
\end{equation}
allows one to express the covariance dynamics (\ref{eq:modified-Kalman-Bucy-covariance-Phit}) in terms of $\mathbf{P}^{-1}(t)$, and then solve for $\mathbf{P}(t)$. Note that (\ref{eq:inverse-covariance-time-derivative}) is easily derived by differentiating the expression $\mathbf{I}_{N_{\mathbf{x}} \times N_{\mathbf{x}}} = \mathbf{P}^{-1}(t) \mathbf{P}(t)$. Substituting (\ref{eq:modified-Kalman-Bucy-covariance-Phit}) into (\ref{eq:inverse-covariance-time-derivative}) gives
\begin{equation} \label{eq:modified-Kalman-Bucy-inverse-covariance-Phi}
    \frac{{\rm d} \mathbf{P}^{-1}(t)}{{\rm d} t} = \frac{1}{\delta} 
    \mathbf{B}  - \mathbf{P}^{-1}(t)
\end{equation}
which has the closed form solution
\begin{equation} \label{eq:modified-Kalman-Bucy-inverse-covariance-Phi-solution}
    \mathbf{P}^{-1}(t) = \exp(-t) \mathbf{P}^{-1}(0) + \frac{1}{\delta} 
    \big ( 1 - \exp(-t) \big ) \mathbf{B} 
\end{equation}
for initial condition $\mathbf{P}^{-1}(0)$. As $\mathbf{B}  \succ 0$, and $\mathbf{P}(0) \succ 0$ implies $\mathbf{P}^{-1}(0) \succ 0$, the right side of (\ref{eq:modified-Kalman-Bucy-inverse-covariance-Phi-solution}) is a weighted sum of positive definite matrices with non-negative weights, and thus also positive definite. Thus, the solution of (\ref{eq:modified-Kalman-Bucy-covariance-Phit}) is
\begin{equation} \label{eq:modified-Kalman-Bucy-covariance-Phi-P-solution}
\mathbf{P}(t) = \delta \bigg [ \big ( 1 - \exp(-t) \big ) \mathbf{B} 
+ \delta \exp(-t) \mathbf{P}^{-1}(0)\bigg ]^{-1} \text{ . }
\end{equation}

Hence, for constant $\mathbf{B} = \mathbf{B}(\mathbf{\bar x})$, the steady state covariance matrix becomes $\mathbf{P}_\infty := \mathbf{P}(t\to \infty) = \delta \mathbf{B}^{-1}$ and (\ref{eq:modified-Kalman-Bucy-covariance-Phit}) provides the desired inverse up to a scaling by $\delta$. Furthermore, upon substituting $\mathbf{P}(t) \approx \delta \mathbf{B}^{-1}$ into (\ref{eq:Kalman-Bucy-mean-evolution-scaled}), we obtain the desired inversion-free reformulation of (\ref{eq:general-inversion-preconditioned-gradient-flow}).

We will later on add random perturbations to particle implementations of (\ref{eq:modified-Kalman-Bucy-covariance-Phit}), which, in the mean-field limit, amounts to replacing (\ref{eq:modified-Kalman-Bucy-covariance-Phit}) by 
\begin{equation} \label{eq:modified-Kalman-Bucy-covariance-noise}
    \frac{{\rm d} \mathbf{P}(t)}{{\rm d}t} = \mathbf{P}(t) - \frac{1}{\delta} \mathbf{P}(t) \mathbf{B} \big ( \mathbf{\bar x}(t) \big ) \mathbf{P}(t) + \beta \delta \mathbf{I}_{N_{\mathbf{x}} \times N_{\mathbf{x}}}
\end{equation}
with parameter $\beta\ge 0$. Provided the stationary $\mathbf{B}$ has an eigenvalue decomposition $\mathbf{B} = \mathbf{V} \mathbf{\Lambda} \mathbf{V}^{\rm T}$, the resulting stationary covariance matrix $\mathbf{P}_\infty$ will be now of the form
\begin{equation}
    \mathbf{P}_\infty = \mathbf{V} \mathbf{D} \mathbf{V}^{\rm T}
\end{equation}
and the entries $d_{ii}$ of the diagonal matrix $\mathbf{D}$ are related to the entries $\lambda_{ii}$ of $\mathbf{\Lambda}$ via
\begin{equation} \label{eq:eigenvalue-noise-relationship}
    d_{ii} = \frac{\delta}{2\lambda_{ii}}\left(
    1+\sqrt{1+4\beta \lambda_{ii}} \right) \text{ . }
\end{equation}
Hence, provided $\beta \ll 1/(4\lambda_{ii})$ for all $i \in \{1,\ldots,N_{\mathbf{x}}\}$, the eigenvalues of $\mathbf{P}_\infty$ are not altered significantly. However, the random perturbations break the well-known subspace property of the EnKF that is present when the standard EnKF is used with fewer particles than the state space dimension $N_{\mathbf{x}}$.

\subsection{Generalized Square Root Formulation Of Covariance}
\label{subsec:covariance-generalized-square-root}
The introduction of a matrix $\mathbf{Y}(t) \in \mathbb{R}^{N_{\mathbf{x}} \times N_{\mathbf{x}}}$ such that $\mathbf{P}(t) = \mathbf{Y}(t) \mathbf{Y}^{\rm T}(t)$ provides a key building block towards EnKF formulations. In other words, $\mathbf{Y}(t)$ provides a generalized square root of the covariance matrix $\mathbf{P}(t)$. Substituting the generalized square root $\mathbf{Y}(t)$ into (\ref{eq:modified-Kalman-Bucy-covariance-Phit}) gives
\begin{equation} \label{eq:modified-Kalman-Bucy-covariance-generalized-square-root-expanded-Phi-P}
    \frac{{\rm d} \mathbf{Y}(t)}{{\rm d} t} \mathbf{Y}^{\rm T}(t) + \mathbf{Y}(t) \frac{{\rm d} \mathbf{Y}^{\rm T}(t)}{{\rm d} t} = \mathbf{Y}(t) \mathbf{Y}^{\rm T}(t) - \frac{1}{\delta} \mathbf{Y}(t) \mathbf{Y}^{\rm T}(t) \mathbf{B} \big ( \mathbf{\bar x}(t) \big ) \mathbf{Y}(t) \mathbf{Y}^{\rm T}(t)
\end{equation}
and (\ref{eq:modified-Kalman-Bucy-covariance-generalized-square-root-expanded-Phi-P}) is satisfied when the following equation holds
\begin{equation} \label{eq:modified-Kalman-Bucy-covariance-generalized-square-root-Phi-P}
    \frac{{\rm d} \mathbf{Y}(t)}{{\rm d} t} = \frac{1}{2} \mathbf{Y}(t) \bigg [ \mathbf{I}_{N_{\mathbf{x}} \times N_{\mathbf{x}}} - \frac{1}{\delta}\mathbf{Y}^{\rm T}(t) \mathbf{B} \big ( \mathbf{\bar x}(t) \big ) \mathbf{Y}(t) \bigg ] \text{ . }
\end{equation}
Note that substituting (\ref{eq:modified-Kalman-Bucy-covariance-generalized-square-root-Phi-P}) into (\ref{eq:modified-Kalman-Bucy-covariance-generalized-square-root-expanded-Phi-P}) verifies that (\ref{eq:modified-Kalman-Bucy-covariance-generalized-square-root-Phi-P}) satisfies (\ref{eq:modified-Kalman-Bucy-covariance-Phit}).

\subsection{Particle Approximation}
\label{subsec:particle-approximation}
Further following the EnKF methodology, we now define an empirical version of $\mathbf{Y}(t)$, which we denote $\mathbf{\hat{Y}}(t) \in \mathbb{R}^{N_{\mathbf{x}} \times K}$, using an ensemble of $K>1$ particles. Let the position of particle $k \in \{1, \dots, K\}$ be denoted by $\mathbf{x}^{(k)}(t) \in \R^{N_{\mathbf{x}}}$ and define the $N_{\mathbf{x}} \times K$ matrix $\mathbf{X}(t)$ of ensemble members as follows:
\begin{equation}
    \mathbf{X}(t) \coloneqq \left( \mathbf{x}^{(1)}(t), \mathbf{x}^{(2)}(t),\ldots, \mathbf{x}^{(K)}(t) \right) \text{ . }
\end{equation}
Let $\hat{m} \big [ \cdot \big ]$ denote the operation of taking the empirical expectation over the particles, so that the empirical mean of the particle positions is
\begin{equation} \label{eq:ensemble_mean}
    \hat{m} \big [ \mathbf{X}(t) \big ] \coloneqq \frac{1}{K} \sum_{k=1}^{K} \mathbf{x}^{(k)}(t) \text{ , }
\end{equation}
which we can rewrite as
\begin{equation}
    \hat{m} \big [ \mathbf{X}(t) \big ] = \mathbf{X}(t) \mathbf{w}_{\rm u} \text{ , }
\end{equation}
where 
\begin{equation}
    \mathbf{w}_{\rm u} := \big ( 1/K, 1/K, \ldots, 1/K \big )^{\rm T} \in \mathbb{R}^K
\end{equation}
is a probability vector of uniform weights $1/K$. We also define the vector of ones
\begin{equation}
    \textbf{1}_K := K \mathbf{w}_{\rm u} = \big ( 1, 1, \ldots, 1 \big )^{\rm T} \in \mathbb{R}^K \text{ . }
\end{equation}
Let the matrix
\begin{equation} \label{eq:ensemble_deviations}
   \mathbf{\tilde{Y}}(t) \coloneqq \bigg ( \mathbf{x}^{(1)}(t) - \hat{m} \big [ \mathbf{x}(t) \big ], \dots, \mathbf{x}^{(K)}(t) - \hat{m} \big [ \mathbf{x}(t) \big ] \bigg ) \in \R^{N_{\mathbf{x}} \times K}
\end{equation}
contain the deviations of the particle positions from the mean position as its columns. Upon introducing the symmetric $K\times K$ projection matrix
\begin{equation}
    \mathbf{\Pi}_K \coloneqq \mathbf{I}_{K\times K} - \mathbf{w}_{\rm u}
    \textbf{1}_{K}^{\rm T} \text{ , }
\end{equation}
which is easily verified to satisfy $\mathbf{\Pi}_K \mathbf{\Pi}_K = \mathbf{\Pi}_K$, (\ref{eq:ensemble_deviations}) can be expressed as 
\begin{equation}
     \mathbf{\tilde{Y}}(t) = \mathbf{X}(t) \mathbf{\Pi}_K \text{ . }
\end{equation}

Define $\mathbf{\hat{Y}}(t) \coloneqq \frac{1}{\sqrt{K}} \mathbf{\tilde{Y}}(t)$. The empirical covariance of the particle positions $\mathbf{\hat{P}}(t)$ is now given by
\begin{equation} \label{eq:empirical-covariance}
    \mathbf{\hat{P}}(t)  \coloneqq
    \frac{1}{K}\mathbf{X}(t) \mathbf{\Pi}_K \mathbf{X}^{\rm T}(t)
    = \frac{1}{K}\mathbf{\tilde{Y}}(t) \mathbf{\tilde{Y}}^{\rm T}(t)
    = \mathbf{\hat{Y}}(t) \mathbf{\hat{Y}}^{\rm T}(t) \text{ , }
\end{equation}
with $\mathbf{\hat{P}}(t) \rightarrow \mathbf{P}(t)$ in the mean-field limit $K \rightarrow \infty$. In other words, $\mathbf{\hat Y}(t)$ provides a particle-based approximation to the generalized square root $\mathbf{Y}(t)$ of $\mathbf{P}(t)$.

We now state evolution equations for the ensemble matrix $\mathbf{X}(t)$ such that the ensemble mean (\ref{eq:ensemble_mean}) and the ensemble deviations (\ref{eq:ensemble_deviations}) satisfy the specified Kalman--Bucy evolution equations (\ref{eq:Kalman-Bucy-mean-evolution-scaled}) and (\ref{eq:modified-Kalman-Bucy-covariance-Phit}) respectively. Thus, we obtain
\begin{equation} \label{eq:ultimate}
    \frac{{\rm d}\mathbf{X}(t)}{{\rm d}t}
    = \frac{1}{2} \mathbf{\tilde Y}(t)  
    - \frac{1}{2\delta K} \mathbf{\tilde Y}(t)
    \mathbf{\tilde Y}^{\rm T}(t)
    \left(
    \mathbf{B} \big ( \mathbf{\bar x}(t) \big ) \mathbf{\tilde Y}(t)
     + 2\nabla_{\bf x} \Phi \big ( \mathbf{\bar x}(t) \big ) \, \textbf{1}_{\rm K}^{\rm T} \right) 
\end{equation}
subject to $\mathbf{\tilde Y}(t) = \mathbf{X}(t)\mathbf{\Pi}_K$ and $\mathbf{\bar x}(t) = \mathbf{X}(t) \mathbf{w}_{\rm u}$. In the next section, we demonstrate that there are implementations of (\ref{eq:ultimate}) which rely only on evaluations of the forward map $\mathcal{G}(\mathbf{x})$, the loss function $\mathcal{D}(\mathbf{y})$ and its known derivatives $\nabla_{\mathbf{y}} \mathcal{D}(\mathbf{y})$ and $\nabla_{\mathbf{y}}^2 \mathcal{D}(\mathbf{y})$, and the regularization functions $\mathcal{R}(\mathbf{x})$ and $\mathcal{T}(\mathbf{y})$ if regularization is present.

\section{Stein's Identity: Particle Based Derivative Approximations}
\label{sec:derivative-approximations}
We now discuss how to use the ensemble of particles to approximate $\mathbf{\hat P}(t) \nabla_{\mathbf{x}} \Phi \big ( \mathbf{\bar{x}}(t) \big )$ and $\mathbf{\hat Y}^{\rm T}(t) \mathbf{B} \big ( \mathbf{\bar x}(t) \big )  \mathbf{\hat Y}(t)$ for different choices of $\mathbf{B} \big ( \mathbf{\bar x}(t) \big )$. We first recall a version of Stein's identity, which we use to estimate unavailable derivatives in the context of (\ref{eq:general-inversion-objective}). 

Let $\mathbf{u} \in \R^{N_{\mathbf{u}}}$ be a random vector and $\mathcal{M}(\mathbf{u}) : \R^{N_{\mathbf{u}}} \mapsto \R^{N_{\mathbf{v}}}$ be a differentiable function. Define the expectation of $\mathcal{M}(\mathbf{u})$ with respect to a probability density function $q(\mathbf{u}) : \R^{N_{\mathbf{u}}} \mapsto \R_{\geq 0}$ as
\begin{equation} \label{eq:expectation-definition}
\mathbb{E}_{q(\mathbf{u})} \big [ \mathcal{M}(\mathbf{u}) \big ] \coloneqq \int \mathcal{M}(\mathbf{u}) \, q(\mathbf{u}) \,{\rm d} \mathbf{u} \text{ . } 
\end{equation}
Now, let $\mathbf{u}$ be multivariate Gaussian distributed $\mathbf{u} \sim \mathcal{N}(\mathbf{\bar u}, \mathbf{\Sigma})$. Denote the corresponding multivariate Gaussian probability density function $p(\mathbf{u})$ by
\begin{equation} \label{eq:Gaussian-probability-density}
p(\mathbf{u}) \coloneqq \frac{1}{\sqrt{(2 \pi)^{N_{\mathbf{u}}} \det \big ( \mathbf{\Sigma} \big )}} \exp \bigg ( - \frac{1}{2} \big ( \mathbf{u} - \mathbf{\bar u} \big )^{\rm T}\mathbf{\Sigma}^{-1} \big ( \mathbf{u} - \mathbf{\bar u} \big ) \bigg ) \text{ . } 
\end{equation}
A version of Stein's identity says that if $\mathbf{u} \sim \mathcal{N}(\mathbf{\bar u}, \mathbf{\Sigma})$, then
\begin{equation} \label{eq:Stein's-identity-gradient}
\mathbf{\Sigma} \,\mathbb{E}_{p(\mathbf{u})} \big [ D_{\mathbf{u}} \mathcal{M}(\mathbf{u} )\big ]^{\rm T} = \mathbb{E}_{p(\mathbf{u})} \big [ \big ( \mathbf{u} - \mathbf{\bar u} \big ) \, \big ( \mathcal{M}(\mathbf{u}) - \mathbb{E}_{p(\mathbf{u})} \big [ \mathcal{M}(\mathbf{u}) \big ] \big )^{\rm T}  \big ] \text{ . }
\end{equation}
The identity (\ref{eq:Stein's-identity-gradient}) is well-known (see \cite{opper-bayesian-online-learning-1999, pidstrigach-ensemble-transform-logistic-2022}), as the proof of (\ref{eq:Stein's-identity-gradient}) via integration by parts is straightforward. Furthermore, the assumption of Gaussianity is justified in our context provided $\delta > 0$ is chosen small enough. 

By the chain rule, the gradient of (\ref{eq:general-inversion-objective}) decomposes as
\begin{equation} \label{eq:general-inversion-objective-gradient-decomposed}
    \nabla_{\mathbf{x}} \Phi(\mathbf{x}) = D_{\mathbf{x}} \mathcal{G} \big ( \mathbf{x} \big )^{\rm T}  \nabla_{\mathbf{y}} \mathcal{D}(\mathbf{y}) + \alpha_{\mathbf{x}} \nabla_{\mathbf{x}} \mathcal{R}(\mathbf{x}) + \alpha_{\mathbf{y}} D_{\mathbf{x}} \mathcal{G} \big ( \mathbf{x} \big )^{\rm T}  \nabla_{\mathbf{y}} \mathcal{T}(\mathbf{y})
\end{equation}
with $\mathbf{y} = \mathcal{G}(\mathbf{x})$. As the Jacobian matrix $D_{\mathbf{x}} \mathcal{G}(\mathbf{x})$ is unavailable, and the gradients $\nabla_{\mathbf{x}} \mathcal{R}(\mathbf{x})$ and $\nabla_{\mathbf{y}} \mathcal{T}(\mathbf{y})$ may also be unavailable, we turn to Monte Carlo approximations based on the version of Stein's identity given in (\ref{eq:Stein's-identity-gradient}). Using the ensemble of $K$ particles, for the forward map $\mathcal{G}(\mathbf{x})$, we define the following quantities
\begin{equation} \label{eq:forward-map-empirical-expectation}
   \mathbb{E}_{p(\mathbf{x},t)} \big [ \mathcal{G} \big ( \mathbf{x} \big ) \big ] \approx \hat{m} \big [ \mathcal{G} \big ( \mathbf{x}(t) \big ) \big ] \coloneqq \frac{1}{K} \sum_{k=1}^{K} \mathcal{G} \big ( \mathbf{x}^{(k)}(t) \big ) \in \R^{N_{\mathbf{y}}}
\end{equation}
\begin{equation} \label{eq:forward-map-gamma-matrix}
    \mathbf{\Gamma}(t) \coloneqq \bigg ( \mathcal{G} \big ( \mathbf{x}^{(1)}(t) \big ) - \hat{m} \big [ \mathcal{G} \big ( \mathbf{x}(t) \big ) \big ], \dots, \mathcal{G} \big ( \mathbf{x}^{(K)}(t) \big ) - \hat{m} \big [ \mathcal{G} \big ( \mathbf{x}(t) \big ) \big ] \bigg ) \in \R^{N_{\mathbf{y}} \times K}  
\end{equation}
where $p(\mathbf{x},t)$ denotes the distribution of a time-varying Gaussian random variable $\mathbf{x}(t)$ with mean $\mathbf{\bar x}(t)$ and covariance $\mathbf{P}(t)$. If the gradient $\nabla_{\mathbf{x}} \mathcal{R}(\mathbf{x})$ of the state space regularization function is unavailable, we also define
\begin{equation} \label{eq:state-space-regularizer-empirical-expectation}
    \mathbb{E}_{p(\mathbf{x},t)} \big [ \mathcal{R} \big ( \mathbf{x} \big ) \big ] \approx \hat{m} \big [ \mathcal{R} \big ( \mathbf{x}(t) \big ) \big ] \coloneqq \frac{1}{K} \sum_{k=1}^{K} \mathcal{R} \big ( \mathbf{x}^{(k)}(t) \big ) \in \R
\end{equation}
\begin{equation} \label{eq:state-space-regularizer-R-vector}
    \mathbf{R}(t) \coloneqq \bigg ( \mathcal{R} \big ( \mathbf{x}^{(1)}(t) \big ) - \hat{m} \big [ \mathcal{R} \big ( \mathbf{x}(t) \big ) \big ], \dots, \mathcal{R} \big ( \mathbf{x}^{(K)}(t) \big ) - \hat{m} \big [ \mathcal{R} \big ( \mathbf{x}(t) \big ) \big ] \bigg ) \in \R^{K}  \text{ . }
\end{equation}
Now, by applying Stein's identity (\ref{eq:Stein's-identity-gradient}) with the choices $\mathbf{u} = \mathbf{x} \sim \mathcal{N}(\mathbf{\bar x}, \mathbf{P})$ and $\mathcal{M} = \mathcal{G}$ (in this case $N_{\mathbf{u}} = N_{\mathbf{x}}$ and $N_{\mathbf{v}} = N_{\mathbf{y}}$), we obtain the Monte Carlo approximation
\begin{equation} \label{eq:Monte-Carlo-Stein-forward-map}
    \mathbf{P}(t) \,\mathbb{E}_{p(\mathbf{x},t)} \big [ D_{\mathbf{x}} 
    \mathcal{G} \big ( \mathbf{x} \big ) \big ]^{\rm T} \approx \frac{1}{K} 
    \mathbf{\tilde{Y}}(t)\,\mathbf{\Gamma}(t)^{\rm T} \in \R^{N_{\mathbf{x}} \times N_{\mathbf{y}}}  \text{ , }
\end{equation}
and by using $\mathcal{M} = \mathcal{R}$ instead (in this case $N_{\mathbf{u}} = N_{\mathbf{x}}$ and $N_{\mathbf{v}} = 1$), we obtain
\begin{equation} \label{eq:Monte-Carlo-Stein-state-space-regularizer}
    \mathbf{P}(t) \,\mathbb{E}_{p(\mathbf{x},t)} \big [ \nabla_{\mathbf{x}} 
    \mathcal{R} \big ( \mathbf{x} \big ) \big ] \approx \frac{1}{K} 
    \mathbf{\tilde{Y}}(t)\, \mathbf{R}(t) \in \R^{N_{\mathbf{x}}}  \text{ . }
\end{equation}
From (\ref{eq:Monte-Carlo-Stein-forward-map}) and (\ref{eq:Monte-Carlo-Stein-state-space-regularizer}), recalling (\ref{eq:empirical-covariance}), we further deduce that
\begin{equation} \label{eq:Monte-Carlo-Y-jacobian}
    \mathbf{\tilde Y}(t) \,\mathbb{E}_{p(\mathbf{x},t)} \big [ D_{\mathbf{x}} \mathcal{G} \big ( \mathbf{x} \big ) \big ] \approx  \mathbf{\Gamma}(t)
\end{equation}
\begin{equation} \label{eq:Monte-Carlo-Y-gradient}
    \mathbf{\tilde Y}(t)^{\rm T} \, \mathbb{E}_{p(\mathbf{x},t)} \big [ \nabla_{\mathbf{x}} \mathcal{R} \big ( \mathbf{x} \big ) \big ] \approx  \mathbf{R}(t)  \text{ . }
\end{equation}
Thus, the terms in (\ref{eq:ultimate}) involving $\nabla_{\mathbf{y}} \mathcal{D}(\mathbf{y})$ and $\nabla_{\mathbf{x}} \mathcal{R}(\mathbf{x})$ can be approximately evaluated using (\ref{eq:Monte-Carlo-Y-jacobian}) and (\ref{eq:Monte-Carlo-Y-gradient}) when the derivatives $D_{\mathbf{x}} \mathcal{G} \big ( \mathbf{x} \big )$ and $\nabla_{\mathbf{x}} \mathcal{R}(\mathbf{x})$ are unavailable. The state space regularization function $\mathcal{T}(\mathbf{y})$ can be treated analogously to the loss function $\mathcal{D}(\mathbf{y})$ if $\nabla_{\mathbf{y}} \mathcal{T}(\mathbf{y})$ is available, and analogously to $\mathcal{R}(\mathbf{x})$ if $\nabla_{\mathbf{y}} \mathcal{T}(\mathbf{y})$ is unavailable.

For second derivatives, the choices for $\mathbf{B} \big ( \mathbf{\bar x}(t) \big )$ we consider are based on the approximate Hessian of (\ref{eq:general-inversion-objective}) given by
\begin{equation} \label{eq:composite-Hessian-B}
    \nabla_{\mathbf{x}}^2 \Phi(\mathbf{x}) \approx D_{\mathbf{x}} \mathcal{G}(\mathbf{x})^{\rm T} \bigg [ \nabla_\mathbf{y}^2
    \mathcal{D}(\mathbf{y}) + \alpha_{\mathbf{y}} \nabla_\mathbf{y}^2 \mathcal{T}(\mathbf{y}) \bigg ] D_{\mathbf{x}} \mathcal{G}(\mathbf{x}) + \alpha_{\mathbf{x}} \nabla_\mathbf{x}^2 \mathcal{R}(\mathbf{x})
\end{equation}
with $\mathbf{y} = \mathcal{G}(\mathbf{x})$. As the Jacobian of the forward map $D_{\mathbf{x}}\mathcal{G}(\mathbf{x})$ is unavailable, instead we will use $\mathbb{E}_{p(\mathbf{x},t)} \left[ D_{\mathbf{x}} \mathcal{G}(\mathbf{x}) \right]$. The choices of $\mathbf{B} \big ( \mathbf{\bar x}(t) \big )$ corresponding to (\ref{eq:composite-Hessian-B}) can be decomposed into the three term sum
\begin{equation} \label{eq:composite-Hessian-B-three-term-sum}
    \mathbf{B} \big ( \mathbf{\bar x}(t) \big ) = \mathbf{B}_{\mathcal{D}} \big ( \mathbf{\bar x}(t) \big ) + \mathbf{B}_{\mathcal{R}} \big ( \mathbf{\bar x}(t) \big ) + \mathbf{B}_{\mathcal{T}} \big ( \mathbf{\bar x}(t) \big )  \text{ . }
\end{equation}
The loss function component $\mathbf{B}_{\mathcal{D}} \big ( \mathbf{\bar x}(t) \big )$ is always given by
\begin{equation} \label{eq:composite-Hessian-B-loss-component}
    \mathbf{B}_{\mathcal{D}} \big ( \mathbf{\bar x}(t) \big ) = \mathbb{E}_{p(\mathbf{x},t)} \left[ D_{\mathbf{x}} \mathcal{G}(\mathbf{x}) \right]^{\rm T} \nabla_\mathbf{y}^2 \mathcal{D} \big ( \mathbf{\bar y}(t) \big ) \, \mathbb{E}_{p(\mathbf{x},t)}\left[ D_{\mathbf{x}} \mathcal{G}(\mathbf{x})\right]  \text{ , }
\end{equation}
where $\mathbf{\bar y}(t) = \mathcal{G} \big ( \mathbf{\bar x}(t) \big )$. If $\nabla_\mathbf{x}^2 \mathcal{R}(\mathbf{x})$ is available, then the state space regularization component is given by
\begin{equation} \label{eq:composite-Hessian-B-state-space-component-differentiable}
    \mathbf{B}_{\mathcal{R}} \big ( \mathbf{\bar x}(t) \big ) = \alpha_{\mathbf{x}} \nabla_\mathbf{x}^2 \mathcal{R} \big ( \mathbf{\bar x}(t) \big )  \text{ , }
\end{equation}
and otherwise $\mathbf{B}_{\mathcal{R}} \big ( \mathbf{\bar x}(t) \big ) = \alpha_{\mathbf{x}} \mathbf{I}_{N_{\mathbf{x}} \times N_{\mathbf{x}}}$. Similarly, if $\nabla_\mathbf{y}^2 \mathcal{T}(\mathbf{y})$ is available, then the observation space regularization component is given by
\begin{equation} \label{eq:composite-Hessian-B-observation-space-component-differentiable}
    \mathbf{B}_{\mathcal{T}} \big ( \mathbf{\bar x}(t) \big ) = \alpha_{\mathbf{y}} \mathbb{E}_{p(\mathbf{x},t)} \left[ D_{\mathbf{x}} \mathcal{G}(\mathbf{x}) \right]^{\rm T} \nabla_\mathbf{y}^2 \mathcal{T} \big ( \mathbf{\bar y}(t) \big ) \, \mathbb{E}_{p(\mathbf{x},t)}\left[ D_{\mathbf{x}} \mathcal{G}(\mathbf{x})\right]  \text{ , }
\end{equation}
and otherwise $\mathbf{B}_{\mathcal{T}} \big ( \mathbf{\bar x}(t) \big ) = \alpha_{\mathbf{y}} \mathbf{I}_{N_{\mathbf{x}} \times N_{\mathbf{x}}}$.

To conclude with a concrete example, the particle-based evolution equations
\begin{equation} \label{eq:ultimate_Stein}
    \frac{{\rm d}\mathbf{X}(t)}{{\rm d}t}
    = \frac{1}{2} \mathbf{\tilde Y}(t)   - \frac{1}{2\delta K}
    \mathbf{\tilde Y}(t) 
    \bigg ( \mathbf{\Gamma}^{\rm T}(t)
    \nabla_\mathbf{y}^2
    \mathcal{D} \big ( \mathbf{\bar y}(t) \big )
    \mathbf{\Gamma}(t)
    + 2 \mathbf{\Gamma}^{\rm T}(t)  \nabla_\mathbf{y} \mathcal{D} \big ( \mathbf{\bar y}(t) \big ) \textbf{1}_{K}^{\rm T} \bigg )
\end{equation}
correspond to the regularization-free case $\alpha_{\mathbf{x}} = \alpha_{\mathbf{y}} = 0$ where $\mathbf{B} \big ( \mathbf{\bar x}(t) \big )$ is given by the single term (\ref{eq:composite-Hessian-B-loss-component}).

\section{EnKSGD Algorithms}
\label{sec:EnKSGD-algorithms}
We are now ready to put everything together into a class of algorithms we call Ensemble Kalman--Stein Gradient Descent (EnKSGD). The only missing ingredient is an appropriate time discretization of the particle-based evolution equations (e.g. (\ref{eq:ultimate_Stein})). Our implementations of EnKSGD use an extension of the well-known ensemble transform Kalman filter (ETKF) \cite{bishop-ETKF-2001}, which we summarize next. In combination with the application of the version of Stein's identity given in (\ref{eq:Stein's-identity-gradient}), we then present a complete EnKSGD algorithm.

\subsection{EnKF Update Step}
\label{subsec:ETKF-update}
To develop an efficient time-stepping method for (\ref{eq:ultimate_Stein}), we first return to the case of a quadratic objective function and a particle implementation of the mean-field equation (\ref{eq:two-scale-EnKBF}). Following our previously derived evolution equations
(\ref{eq:ultimate}), and by noting that $\mathbf{B} \big ( \mathbf{\bar x}(t) \big ) = \mathbf{G}^{\rm T}\mathbf{G}$ for the quadratic objective function (\ref{eq:quadratic-objective-function-phi}), we obtain the continuous-time formulation
\begin{equation}
    \frac{{\rm d}\mathbf{X}(t)}{{\rm d}t} =
    \frac{1}{2} \mathbf{\tilde Y}(t)  
    - \frac{1}{2\delta K}\mathbf{\tilde Y}(t)
    \mathbf{\tilde Y}^{\rm T}(t) \left(
    \mathbf{G}^{\rm T} \mathbf{G} \mathbf{\tilde Y}(t)
    + 2\mathbf{G}^{\rm T}(\mathbf{G}\mathbf{\bar x}(t) - \mathbf{y}_{\rm obs}) \textbf{1}_{K}^{\rm T} 
    \right).
\end{equation}

An efficient time-stepping method is provided by the 
ensemble transform Kalman filter (ETKF)
applied to the continuous-time Kalman--Bucy filter equations \cite{amezcua-ETKBF-2014}. More specifically, given an initial ensemble $\mathbf{X}_0 \in \mathbb{R}^{N_{\mathbf{x}} \times K}$ and sequence of step sizes $\Delta t_n >0$, the ETKF proceeds as follows. For each $n\ge 0$, the ensemble deviations $\mathbf{\tilde Y}_n = \mathbf{X}_n \mathbf{\Pi}_K$ are updated from time $t_n$ to 
$t_{n+1} = t_n + \Delta t_n$ as follows:
\begin{equation} \label{eq:ETKF_deviations}
    \mathbf{\tilde Y}_{n+1} = \exp \bigg ( \frac{\Delta t_n}{2} \bigg ) \mathbf{\tilde Y}_n \mathbf{T}_n^{1/2} \text{ . }
\end{equation}
Here, the $K\times K$ matrix $\mathbf{T}_n$ is defined as
\begin{equation}
    \mathbf{T}_n \coloneqq \left( \mathbf{I}_{K\times K} + \frac{\Delta t_n}{\delta K}
    \mathbf{\tilde Y}_n^{\rm T} \mathbf{G}^{\rm T} \mathbf{G}
    \mathbf{\tilde Y}_n \right)^{-1}.
\end{equation}
The update for the ensemble mean $\mathbf{\bar x}_n = \mathbf{X}_n \mathbf{w}_{\rm u}$ is given by 
\begin{subequations}
\begin{align}
    \mathbf{\bar x}_{n+1} &= 
    \mathbf{\bf X}_n
    \left( \mathbf{w}_{\rm u} - \frac{\Delta t_n}{\delta K} \mathbf{\Pi}_K \mathbf{T}_n
    \mathbf{\tilde Y}_n^{\rm T} \mathbf{G}^{\rm T}
    \left( \mathbf{G}\mathbf{\bar x}_n - \mathbf{y}_{\rm obs}\right)\right) \\ \label{eq:ETKF_mean}
    &=
    \mathbf{\bf X}_n
    \left( \mathbf{w}_{\rm u} - \frac{\Delta t_n}{\delta K} \mathbf{\Pi}_K \mathbf{T}_n
    \mathbf{\tilde Y}_n^{\rm T} \nabla_{\bf x} \Phi(\mathbf{\bar x}_n) \right).
\end{align}
\end{subequations}

We observe that the only ingredient that requires adaptation to an objective function with the structure of (\ref{eq:general-inversion-objective}) is the definition of the transform matrix $\mathbf{T}_n$. For general $\mathbf{B} \big ( \mathbf{\bar{x}}(t) \big )$, a natural choice is
\begin{equation} \label{eq:ETKF_transform_matrix}
    \mathbf{T}_n = \left( \mathbf{I}_{K\times K} + \frac{\Delta t_n}{\delta K} \mathbf{\tilde Y}_n^{\rm T} 
    \mathbf{B}(\mathbf{\bar x}_n)
    \mathbf{\tilde Y}_n \right)^{-1}.
\end{equation}
Using Monte Carlo approximations of the form (\ref{eq:Monte-Carlo-Y-jacobian}) and (\ref{eq:Monte-Carlo-Y-gradient}), we obtain the proposed EnKSGD class of algorithms, which we summarize in the following subsection.

\subsection{EnKSGD Update Step}
\label{subsec:EnKIPGD-update-steps}
Combining the particle-based evolution equations (\ref{eq:ultimate_Stein}) with the time discretization (\ref{eq:ETKF_mean}) leads to the derivative-free update
\begin{equation} \label{eq:EnKIPGD_mean}
     \mathbf{\bar x}_{n+1} = \mathbf{\bf X}_n
     \left( \mathbf{w}_{\rm u} - \frac{\Delta t_n}{\delta K} \mathbf{\Pi}_K \mathbf{T}_n \mathbf{\Gamma}_n^{\rm T}
     \nabla_\mathbf{y} \mathcal{D}(\mathbf{\bar y}_n) \right)
\end{equation}
of the ensemble mean $\mathbf{\bar x}_n$, where $\mathbf{\bar y}_n = \mathcal{G}(\mathbf{\bar x}_n)$ and
\begin{equation} \label{eq:forward-map-gamma-matrix-discretized}
     \mathbf{\Gamma}_n \coloneqq \bigg ( \mathcal{G} \big ( \mathbf{x}^{(1)}_n \big ) - \hat{m} \big [ \mathcal{G} \big ( \mathbf{x}_n \big ) \big ], \dots, \mathcal{G} \big ( \mathbf{x}^{(K)}_n \big ) - \hat{m} \big [ \mathcal{G} \big ( \mathbf{x}_n \big ) \big ] \bigg )  \text{ . }
\end{equation}
In place of (\ref{eq:ETKF_transform_matrix}), the ensemble transform matrix $\mathbf{T}_n$ is now given by
\begin{equation} \label{eq:EnKSGD_transform_matrix}
     \mathbf{T}_n = \left( \mathbf{I}_{K\times K} + \frac{\Delta t_n}{\delta K} \mathbf{\Gamma}_n^{\rm T} 
     \nabla_\mathbf{y}^2
     \mathcal{D}(\mathbf{\bar y}_n)
     \mathbf{\Gamma}_n 
     \right)^{-1}  . 
\end{equation}
Combining the update (\ref{eq:ETKF_deviations}) for $\mathbf{\tilde Y}_{n+1}$ and the update (\ref{eq:EnKIPGD_mean}) for $\mathbf{\bar x}_{n+1}$ into a single update for $\mathbf{X}_{n+1}$ produces
\begin{equation}
\mathbf{X}_{n+1} = \mathbf{\tilde Y}_{n+1} + \mathbf{\bar x}_{n+1} \textbf{1}_K^{\rm T}  \text{ . }
\end{equation}
The EnKSGD updating step for (\ref{eq:ultimate_Stein}) can therefore be concisely summarized as
\small
\begin{equation} \label{eq:EnKSGD_step}
     \mathbf{X}_{n+1} = \mathbf{X}_{n}\bigg ( \exp \bigg ( \frac{\Delta t_{n}}{2} \bigg ) \mathbf{\Pi}_K {\bf T}_n^{1/2} + \bigg (
     \mathbf{w}_{\rm u}
     - \frac{\Delta t_n}{\delta K} \mathbf{\Pi}_K {\bf T}_n 
     \mathbf{\Gamma}_n^{\rm T}
     \nabla_\mathbf{y} \mathcal{D}(\mathbf{\bar y}_n) \bigg )
     \textbf{1}_K^{\rm T} \bigg )  \text{ . }
\end{equation}
\normalsize
The ensemble transform matrix $\mathbf{T}_n$ is provided by either (\ref{eq:ETKF_transform_matrix}) or (\ref{eq:EnKSGD_transform_matrix}) respectively. 

It should be noted that the step size $\Delta t_n$ has to be chosen such that $\mathbf{T}_n$ remains positive definite. This restriction does not apply to (\ref{eq:EnKSGD_transform_matrix}) but may apply to (\ref{eq:ETKF_transform_matrix}) in case $\mathbf{B}(\mathbf{\bar x}_n)$ is not positive definite.

We finally need to break the subspace property in the case where $K < N_{\mathbf{x}}$. One possibility is to add random perturbations, 
as already discussed in (\ref{eq:modified-Kalman-Bucy-covariance-noise}), leading to the following modification of (\ref{eq:EnKSGD_step}):
\small
\begin{equation*} \label{eq:EnKSGD_step_random}
     \mathbf{X}_{n+1} = \mathbf{X}_{n} \bigg ( \exp \bigg ( \frac{\Delta t_{n}}{2} \bigg ) \mathbf{\Pi}_K {\bf T}_n^{1/2} + \bigg (
     \mathbf{w}_{\rm u}
     - \frac{\Delta t_n}{\delta K} \mathbf{\Pi}_K {\bf T}_n 
     \mathbf{\Gamma}_n^{\rm T}
     \nabla_\mathbf{y} \mathcal{D}(\mathbf{\bar y}_n) \bigg )
     \textbf{1}_K^{\rm T} \bigg ) + \sqrt{\beta \delta \Delta t_{n}}\,
     \mathbf{\Xi}_n \mathbf{\Pi}_K
\end{equation*}
\normalsize
where $\mathbf{\Xi}_n \in \R^{N_{\mathbf{x}} \times K}$ contains $K$ independent perturbations drawn from the $N_{\mathbf{x}}$ dimensional standard Gaussian distribution $\mathcal{N}(\mathbf{0}, \mathbf{I}_{N_{\mathbf{x}} \times N_{\mathbf{x}}})$ and $\beta \ge 0$ is a tuneable parameter. The ensemble mean is preserved via multiplying the perturbations $\mathbf{\Xi}_n$ by $\mathbf{\Pi}_K$. Adding the random perturbations breaks the subspace property of EnKF-type algorithms but also violates the principle of affine invariance. However, the impact is minimal provided $\beta >0$ is chosen small enough (see equation (\ref{eq:eigenvalue-noise-relationship}) in subsection~\ref{subsec:two-time-scale-approach}).

\subsection{EnKSGD Minimization Procedure}
\label{subsec:EnKSGD-minimization-procedure}

\begin{algorithm}[H]
\caption{EnKSGD Minimization Procedure}
\label{alg:EnKSGD-minimization}
\begin{algorithmic}[1]
\Procedure{EnKSGD-Min}{$\mathbf{\bar{x}}_{0}, \mathbf{\tilde{Y}}_{0}$, $\beta$, $\delta$, $\mu_{ls}$, $c_{ls}$, $\tau_{ls}$, $\ell_{\rm max}$, $\gamma_{ub}$, $\gamma_{lb}$, $N_{\rm max}$}
\smallskip
\State Set the initial particle locations $\mathbf{X}_{0} \gets \mathbf{\tilde{Y}}_{0} \mathbf{\Pi}_K + \mathbf{\bar{x}}_{0} \textbf{1}_{K}^{\rm T}$
\smallskip
\For{$n = 0, \dots, N_{\rm max}$}
\smallskip
\State Calculate the forward map deviations matrix $\mathbf{\Gamma}_{n}$ using (\ref{eq:forward-map-gamma-matrix-discretized}) 
\smallskip
\State Calculate $\mathbf{q}_{n} \coloneqq \mathbf{\tilde Y}_n^{\rm T} \nabla_{\bf x} \Phi(\mathbf{\bar x}_n)$ using (\ref{eq:Monte-Carlo-Y-jacobian}) and (\ref{eq:Monte-Carlo-Y-gradient}) as needed
\smallskip
\State Set $\Delta t_{n}^{\prime} \gets \mu_{\rm ls}$, $\mathcal{I}_{ls} \gets {\rm False}$, $\ell \gets 0$ 
\smallskip
\While{$\mathcal{I}_{ls}$ \textbf{is} ${\rm False}$}
\smallskip
\If{$\ell \geq \ell_{\rm max}$}
\smallskip
\State Set $\Delta t_{n} \gets 0$, $\mathbf{\bar{x}}_{n+1} \gets \mathbf{\bar{x}}_{n}$, $\mathbf{T}_{n} \gets \mathbf{I}_{K \times K}$
\smallskip
\State \textbf{break}
\smallskip
\EndIf
\State Calculate the proposed transform matrix $\mathbf{T}_{n}^{\prime}$ using (\ref{eq:ETKF_transform_matrix}) or (\ref{eq:EnKSGD_transform_matrix}) \label{EnKSGD-procedure-proposed-transform-matrix}
\smallskip
\State Calculate $\mathbf{r}_{n}^{\prime} \coloneqq \frac{\Delta t_{n}^{\prime}}{\delta K} \mathbf{T}_{n}^{\prime} \mathbf{q}_{n}$
\smallskip
\State Propose the ensemble mean update $\mathbf{\bar{x}}_{n+1}^{\prime} \gets \mathbf{\bar{x}}_{n} - \mathbf{\tilde Y}_n \mathbf{r}_n^{\prime}$
\smallskip
\If{$\Phi(\mathbf{\bar{x}}_{n+1}^{\prime}) \leq \Phi(\mathbf{\bar{x}}_{n}) - c_{ls} \mathbf{q}_{n}^{\rm T} \mathbf{r}_{n}^{\prime}$} \label{EnKSGD-procedure-check-approx-suff-dec}
\smallskip
\State Set $\mathcal{I}_{ls} \gets {\rm True}$, $\Delta t_{n} \gets \Delta t_{n}^{\prime}$, $\mathbf{\bar{x}}_{n+1} \gets \mathbf{\bar{x}}_{n+1}^{\prime}$, $\mathbf{T}_{n} \gets \mathbf{T}_{n}^{\prime}$
\smallskip
\Else
\smallskip
\State Set $\Delta t_{n}^{\prime} \gets \tau_{ls} \Delta t_{n}^{\prime}$, $\ell \gets \ell + 1$
\smallskip
\EndIf
\EndWhile
\State Update the deviations $\mathbf{\tilde Y}_{n+1} \gets \exp \big ( \frac{\Delta t_{n}}{2} \big ) \mathbf{\tilde Y}_{n} \mathbf{T}_n^{1/2} + \sqrt{\beta \delta \Delta t_{n}} \mathbf{\Xi}_n$ \label{EnKSGD-procedure-deviations-update-only}
\smallskip
\For{$k = 1, \dots, K$} \label{EnKSGD-procedure-deviations-bound-check-start}
\smallskip
\If{$\frac{1}{N_{\mathbf{x}}} \norm{\mathbf{\tilde Y}_{n+1}[:,k]}_2 > \gamma_{ub}$}
\State Set $\mathbf{\tilde Y}_{n+1}[:,k] \gets \gamma_{ub} \mathbf{\tilde Y}_{n+1}[:,k] \bigg / \norm{\mathbf{\tilde Y}_{n+1}[:,k]}_2$
\EndIf
\If{$\frac{1}{N_{\mathbf{x}}} \norm{\mathbf{\tilde Y}_{n+1}[:,k]}_2 < \gamma_{lb}$}
\State Set $\mathbf{\tilde Y}_{n+1}[:,k] \gets \gamma_{lb} \mathbf{\tilde Y}_{n+1}[:,k] \bigg / \norm{\mathbf{\tilde Y}_{n+1}[:,k]}_2$
\EndIf
\EndFor \label{EnKSGD-procedure-deviations-bound-check-stop}
\State Update the particle locations $\mathbf{X}_{n+1} \gets \mathbf{\tilde{Y}}_{n+1} \mathbf{\Pi}_K + \mathbf{\bar{x}}_{n+1} \textbf{1}_{K}^{\rm T}$ \label{update-particle-locations}
\smallskip
\EndFor
\State \textbf{return} $\Phi(\mathbf{\bar{x}}_{n+1})$
\EndProcedure
\end{algorithmic}
\end{algorithm}

Our EnKSGD minimization procedure proceeds as shown in Algorithm~\ref{alg:EnKSGD-minimization}. As input, the procedure takes an initial mean $\mathbf{\bar{x}}_{0}$ and deviations matrix $\mathbf{\tilde{Y}}_{0}$, the tuneable perturbation parameter $\beta \geq 0$ and scale parameter $\delta > 0$, the backtracking line search initial step size $\mu_{ls} > 0$, decrease factor $0 < c_{ls} < 1$, backtracking factor $0 < \tau_{ls} < 1$, and maximum number of backtracks $\ell_{\rm max}$, the deviation bounds $0 \leq \gamma_{lb} < \gamma_{ub}$, and a maximum number of iterations $N_{\rm max}$. The lower bound $\gamma_{lb}$ and upper bound $\gamma_{ub}$ control a safety mechanism in lines \ref{EnKSGD-procedure-deviations-bound-check-start} - \ref{EnKSGD-procedure-deviations-bound-check-stop} of Algorithm~\ref{alg:EnKSGD-minimization} that can be used to clip each column of $\mathbf{\tilde Y}_{n+1}$ if the column becomes dangerously large or small. To preserve the ensemble mean, it is essential that $\mathbf{\tilde{Y}}_{n+1} \mathbf{\Pi}_K$ is used in line~\ref{update-particle-locations} instead of $\mathbf{\tilde{Y}}_{n+1}$.

As a convention, in Algorithm~\ref{alg:EnKSGD-minimization}, quantities marked with ${}^{\prime}$ represent proposed quantities which 
may or may not be accepted. To illustrate, the \textbf{if} statement on line \ref{EnKSGD-procedure-check-approx-suff-dec} of Algorithm~\ref{alg:EnKSGD-minimization} uses Stein's identity to approximately test if the sufficient decrease (i.e. Armijo) condition 
\begin{equation} \label{eq:Phi-sufficient-decrease-condition}
    \Phi(\mathbf{\bar{x}}_{n+1}^{\prime}) \leq \Phi(\mathbf{\bar{x}}_{n}) - c_{ls} \nabla_{\mathbf{x}} \Phi(\mathbf{\bar x}_n)^{\rm T} \mathbf{\tilde Y}_n \mathbf{r}_{n}^{\prime}
\end{equation}
is satisfied for a proposed step size $\Delta t_{n}^{\prime}$.

\section{Numerical Experiments}
\label{sec:numerical-experiments}
In this section, we use Julia \cite{bezanson-julia-2017} to study the empirical performance of Algorithm~\ref{alg:EnKSGD-minimization}. In all experiments, we compare EnKSGD to a more standard EnKF-type approach that replaces line~\ref{EnKSGD-procedure-deviations-update-only} in Algorithm~\ref{alg:EnKSGD-minimization} with
\begin{equation} \label{eq:EnKF-procedure-deviations-update-only}
    \mathbf{\tilde Y}_{n+1} \gets \mathbf{\tilde Y}_{n} \mathbf{T}_n^{1/2} + \sqrt{\beta \delta \Delta t_{n}} \mathbf{\Xi}_n
\end{equation}
and thus, instead of (\ref{eq:modified-Kalman-Bucy-covariance-noise}), corresponds to the covariance dynamics given by
\begin{equation} \label{eq:EnKF-procedure-modified-Kalman-Bucy-covariance-noise}
    \frac{{\rm d} \mathbf{P}(t)}{{\rm d}t} = - \frac{1}{\delta} \mathbf{P}(t) \mathbf{B} \big ( \mathbf{\bar x}(t) \big ) \mathbf{P}(t) + \beta \delta \mathbf{I}_{N_{\mathbf{x}} \times N_{\mathbf{x}}} \text{ . }
\end{equation}

We calculate $\mathbf{T}_{n}^{\prime}$ in line~\ref{EnKSGD-procedure-proposed-transform-matrix} of Algorithm~\ref{alg:EnKSGD-minimization} using the second derivative approximations given in equations (\ref{eq:composite-Hessian-B}) to (\ref{eq:composite-Hessian-B-observation-space-component-differentiable}) from Section~\ref{sec:derivative-approximations}. For the problems we consider, $\mathbf{T}_{n}$ is positive definite, and thus has the decomposition $(\mathbf{T}_{n})^{-1} = \mathbf{U} \mathbf{S} \mathbf{U}^{\rm T}$ where $\mathbf{U}$ is an orthonormal matrix and $\mathbf{S}$ is a diagonal matrix, yielding $\mathbf{T}_{n} = \mathbf{U} \mathbf{S}^{-1} \mathbf{U}^{\rm T}$ and $\mathbf{T}_{n}^{1/2} = \mathbf{U} \mathbf{S}^{-1/2} \mathbf{U}^{\rm T}$. Our implementations of Algorithm~\ref{alg:EnKSGD-minimization} compute $\mathbf{U}$ and $\mathbf{S}$ using the singular value decomposition (SVD) of $(\mathbf{T}_{n})^{-1}$, and add the small positive constant $10^{-7}$ to the diagonal entries of $\mathbf{S}$ to avoid numerical instability. 

In all experiments, the line search parameters are set to $\mu_{ls} = 1$, $c_{ls} = 10^{-4}$, $\tau_{ls} = 0.1$, and $\ell_{\rm max} = 15$, the deviation bounds to $\gamma_{lb} = 10^{-4}$ and $\gamma_{ub} = 10^{4}$, and the initial deviations matrix $\mathbf{\tilde{Y}}_{0}$
contains $K$ independent perturbations drawn from the Gaussian distribution $\mathcal{N}(\mathbf{0}, \sigma_{0}^2 \mathbf{I}_{N_{\mathbf{x}} \times N_{\mathbf{x}}})$ with $\sigma_{0} = 10^{-2}$.

\subsection{Ill-Conditioned Linear Least Squares}
\label{subsec:ill-conditioned-linear-least-squares}
\begin{figure}[H]
  \begin{center}
    {\def\arraystretch{1.5}\tabcolsep=3pt
    \begin{tabular}{| c | c |}
    \hline
    \includegraphics[width=5.5cm]{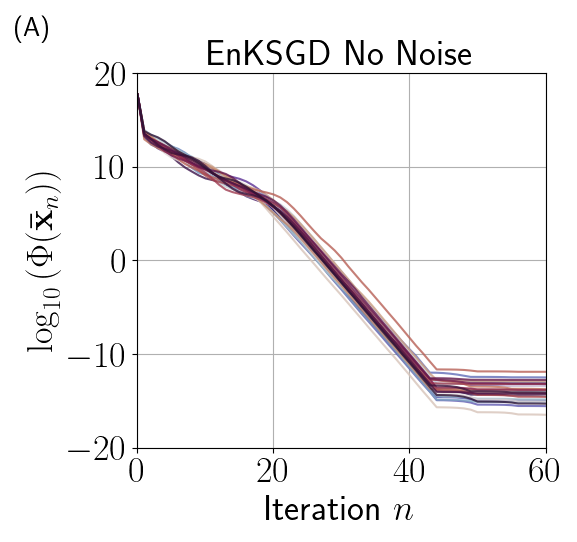} &
    \includegraphics[width=5.5cm]{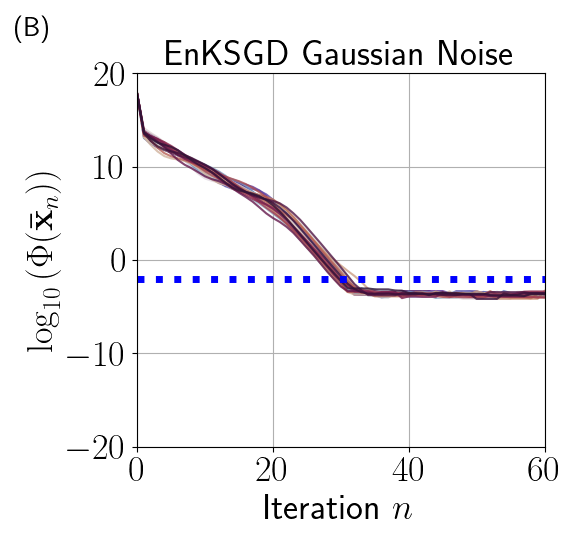} \\ 
    \hline
    \includegraphics[width=5.5cm]{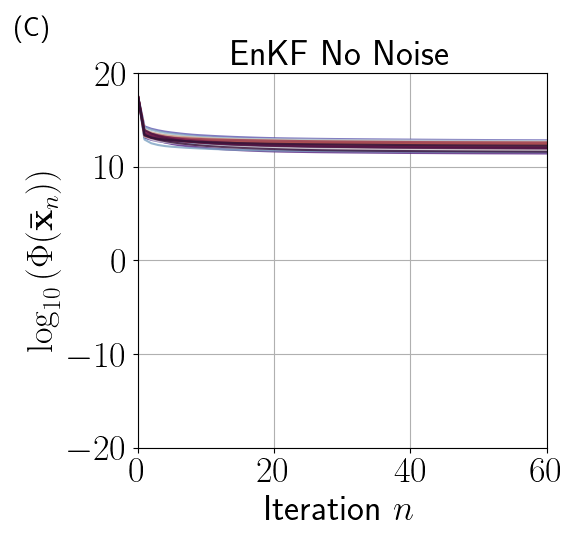} &
    \includegraphics[width=5.5cm]{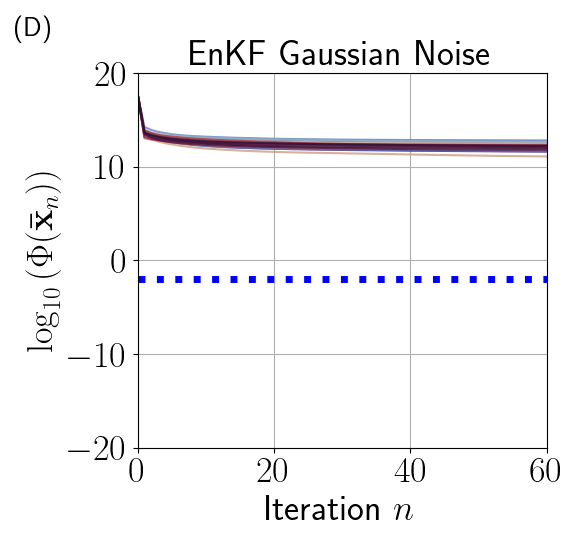} \\
    \hline
    \includegraphics[width=5.5cm]{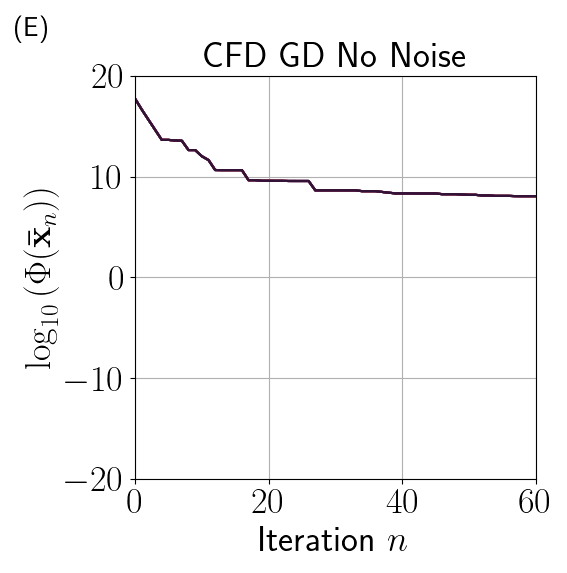} &
    \includegraphics[width=5.5cm]{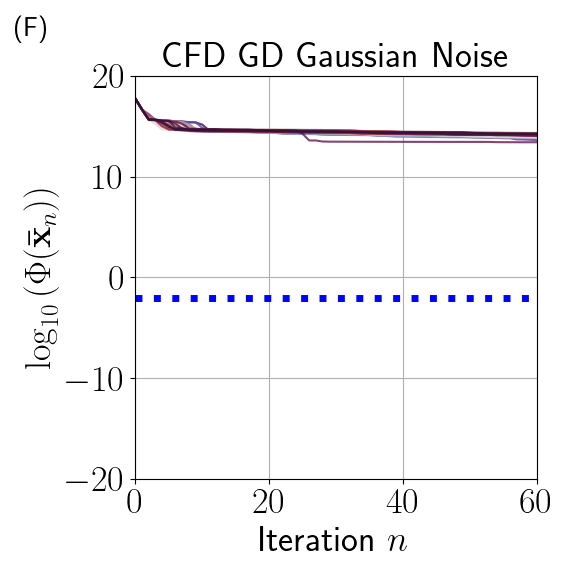} \\
    \hline
    \end{tabular}
    \caption{\textbf{Linear Least Squares}: Results of $30$ independent runs of EnKSGD, EnKF, and central finite difference gradient descent for an ill-conditioned linear least squares problem with and without IID Gaussian noise added to evaluations of $\mathcal{G}(\mathbf{x})$. The blue dashed lines in (B), (D), and (F) mark the standard deviation $\sigma = 10^{-2}$ of the noise added to $\mathcal{G}(\mathbf{x})$. The average total number of $\mathcal{G}(\mathbf{x})$ evaluations for EnKSGD in (A) and (B) is $1261$ and $1421$ respectively, for EnKF in both (C) and (D) is $1261$, and for central finite difference gradient descent in (E) and (F) is $1885$ and $2067$. }
    \label{fig:linear-least-squares-all}
    }
  \end{center}
\end{figure}

The first problem we consider is an ill-conditioned linear least squares problem constructed by choosing $\alpha_{\mathbf{x}} = \alpha_{\mathbf{y}} = 0$, $\mathcal{D}(\mathbf{y}) = \frac{1}{2} \norm{\mathbf{y}}_2^2$, and $\mathcal{G}(\mathbf{x}) = \mathbf{G} \mathbf{x}$. We study the case where the matrix $\mathbf{G}$ is diagonal with entries $g_{ii} = 10^{-2 + 0.5(i-1)}$ for $i \in \{1, \dots, N_{\mathbf{x}}\}$ and $N_{\mathbf{x}} = N_{\mathbf{y}} = 13$. We also study adding noise $\eta(\mathbf{x})$ distributed according to (\ref{eq:stochastic-noise-example}) with $\sigma = 10^{-2}$ to $\mathcal{G}(\mathbf{x})$. The number of particles $K = 20$, the parameters $\beta = 10^{-8}$ and $\delta = 1$, and the initial mean $\mathbf{\bar{x}}_{0} = 10^{5} \cdot \textbf{1}_{N_{\mathbf{x}}}$. Figure~\ref{fig:linear-least-squares-all} shows the progress of EnKSGD, the EnKF-type approach using (\ref{eq:EnKF-procedure-deviations-update-only}), and standard gradient descent (GD) using a central finite difference (CFD) gradient estimate. For all components $i \in \{ 1, \dots, N_{\mathbf{x}} \}$, the CFD gradient estimate is
\begin{equation} \label{eq:CFD-gradient-def}
    \big [ \nabla_{\rm CFD} \Phi(\mathbf{\bar x}_n, h) \big ]_{i} \coloneqq \frac{\Phi(\mathbf{\bar x}_n + h \mathbf{e}_{i}) - \Phi(\mathbf{\bar x}_n - h \mathbf{e}_{i} )}{2 h}
\end{equation}
where $\mathbf{e}_{i}$ denotes a vector with $1$ in component $i$ and $0$ in all other components. For the CFD GD approach, the stencil size is set to $h = 10^{-4}$, and the analogue of the approximate sufficient decrease condition test (\ref{eq:Phi-sufficient-decrease-condition}) is the following test
\begin{equation} \label{eq:CFD-GD-Phi-sufficient-decrease-condition}
    \Phi(\mathbf{\bar{x}}_{n+1}^{\prime}) \leq \Phi(\mathbf{\bar{x}}_{n}) - c_{ls} \Delta t_{n}^{\prime} \nabla_{\rm CFD} \Phi(\mathbf{\bar x}_n, h)^{\rm T} \nabla_{\rm CFD} \Phi(\mathbf{\bar x}_n, h) \text{ . }
\end{equation}

\subsection{Nonlinear Least Squares}
\label{subsec:nonlinear-least-squares}
Next, we consider a variety of nonlinear least squares problems constructed by choosing $\alpha_{\mathbf{x}} = \alpha_{\mathbf{y}} = 0$, $\mathcal{D}(\mathbf{y}) = \frac{1}{2} \norm{\mathbf{y} - \mathbf{y}_{\rm obs}}_2^2$, and with $\mathbf{y}_{\rm obs}$ and $\mathcal{G}(\mathbf{x})$ being problem specific. All problems are taken from \cite{orban-siqueira-nlsproblems-2021}, and the initial means $\mathbf{\bar{x}}_{0}$ are problem specific. The number of particles $K = 8$, $\beta = 10^{-8}$, and $\delta = 10^{-3}$. Table~\ref{tab:nls-type-problems} compares EnKSGD vs. the EnKF-type approach for a fixed budget of $500$ forward map $\mathcal{G}(\mathbf{x})$ evaluations (the number of iterations $n$ can vary).

\begin{table}[h]
    \begin{center}
    {\def\arraystretch{1.5}\tabcolsep=2pt
    \begin{tabular}{| c | c | c | c | c | c | c | c | c |}
        \hline
        \multicolumn{2}{| c |}{\textbf{Characteristics}} & \multicolumn{3}{| c |}{\textbf{EnKSGD} $\log_{10} \big ( \Phi(\mathbf{\bar{x}}) \big )$} & \multicolumn{3}{| c |}{\textbf{EnKF} $\log_{10} \big ( \Phi(\mathbf{\bar{x}}) \big )$} \\
        \hline
        \textit{Problem} & $N_{\mathbf{x}}$ & Mean & Median & Var. & Mean & Median & Var. \\
        \hline
        nls\_rosenbrock & $2$ & $-2.1$E$+01$ & $-2.0$E$+01$ & $2.6$E$+00$ & $+1.4$E$-01$ & $+1.3$E$-01$ & $9.2$E$-04$ \\
        \hline        
        hs25 & $3$ & $+7.8$E$-01$ & $+1.2$E$+00$ & $1.8$E$+00$ & $+1.2$E$+00$ & $+1.2$E$+00$ & $4.6$E$-31$ \\
        \hline
        mgh11 & $3$ & $+4.7$E$-01$ & $+4.8$E$-01$ & $3.8$E$-03$ & $+6.7$E$-01$ & $+6.7$E$-01$ & $6.3$E$-06$ \\
        \hline
        mgh18 & $6$ & $-2.2$E$+00$ & $-2.3$E$+00$ & $6.9$E$-01$ & $-7.9$E$-01$ & $-8.1$E$-01$ & $2.0$E$-03$ \\
        \hline
        tp294 & $6$ & $-1.0$E$+01$ & $-1.2$E$+01$ & $2.0$E$+01$ & $+6.9$E$-01$ & $+6.4$E$-01$ & $3.3$E$-02$ \\
        \hline
        mgh19 & $11$ & $-6.3$E$-01$ & $-6.6$E$-01$ & $7.0$E$-02$ & $-1.1$E$-01$ & $-1.2$E$-01$ & $3.8$E$-02$ \\
        \hline
        tp296 & $16$ & $+3.0$E$+00$ & $+3.0$E$+00$ & $3.2$E$-02$ & $+3.0$E$+00$ & $+3.0$E$+00$ & $4.7$E$-02$ \\
        \hline
        mgh22 & $20$ & $+2.3$E$+00$ & $+2.3$E$+00$ & $3.1$E$-02$ & $+2.3$E$+00$ & $+2.3$E$+00$ & $1.8$E$-02$ \\
        \hline
        tp297 & $30$ & $+3.8$E$+00$ & $+3.8$E$+00$ & $8.8$E$-03$ & $+3.9$E$+00$ & $+3.9$E$+00$ & $1.2$E$-02$ \\
        \hline
        tp304 & $50$ & $+4.0$E$-01$ & $+3.3$E$-01$ & $2.8$E$-01$ & $+3.5$E$+00$ & $+3.5$E$+00$ & $1.2$E$-01$ \\
        \hline
        tp305 & $100$ & $+1.4$E$+00$ & $+1.2$E$+00$ & $8.0$E$-01$ & $+4.7$E$+00$ & $+4.7$E$+00$ & $1.3$E$-01$ \\
        \hline
    \end{tabular}
    \caption{\textbf{Nonlinear Least Squares}: Performance on $11$ nonlinear least squares problems of varying dimension after $500$ evaluations of $\mathcal{G}(\mathbf{x})$ for both EnKSGD and the EnKF-type approach using (\ref{eq:EnKF-procedure-deviations-update-only}). Statistics are calculated with samples of $30$ independent runs, and ``Var.'' denotes variance. }
    \label{tab:nls-type-problems}
    }
    \end{center}
\end{table}

\subsection{Poisson Regression}
\label{subsec:poisson-regression}
We now investigate a maximum likelihood estimation (MLE) problem, where the goal is to minimize the negative log likelihood (NLL) of $N_{\mathbf{y}}$ observed data points given the model parameters $\mathbf{x}$. The forward map $\mathcal{G}(\mathbf{x})$ outputs a vector with $N_{\mathbf{y}}$ components, each representing the probability of an individual data point given $\mathbf{x}$. Poisson regression is a type of generalized linear model (GLM) that models count data using the Poisson distribution; the forward map is
\begin{equation} \label{eq:poisson-regression-forward-map}
    \big [ \mathcal{G}(\mathbf{x}) \big ]_{i} = \frac{\exp(b_{i} \mathbf{x}^{\rm T} \mathbf{a}_{i}) \exp(-\exp(\mathbf{x}^{\rm T} \mathbf{a}_{i}))}{b_i !} \text{ , } \quad \forall i \in \{ 1, \dots, N_{\mathbf{y}} \}
\end{equation}
where $b_{i} \in \mathbb{N}_{0}$ is an observed count (e.g. $0$, $1$, $2$, etc.) and $\mathbf{a}_{i} \in \mathbb{R}^{N_{\mathbf{x}}}$. To compute the NLL, $\mathcal{D}(\mathbf{y}) = - \textbf{1}_{N_{\mathbf{y}}}^{\rm T} \log(\mathbf{y})$ where $\log(\mathbf{y}) \in \mathbb{R}^{N_{\mathbf{y}}}$ denotes the natural logarithm applied component-wise to $\mathbf{y}$. This problem is not regularized, so $\alpha_{\mathbf{x}} = \alpha_{\mathbf{y}} = 0$.

We simulated a Poisson regression dataset with $N_{\mathbf{x}} = 41$ features by first drawing $N_{\mathbf{y}} = 189$ samples $\mathbf{a}_{i} \sim \mathcal{N}(\mathbf{0}, \mathbf{\Sigma})$ with $\mathbf{\Sigma}$ chosen as the ill-conditioned diagonal matrix with diagonal entries given by $10^{-5 + 0.1 (m-1)}, m \in \{ 1, \dots, N_{\mathbf{x}} \}$. Second, we drew $b_i \sim {\rm Pois}(\exp(\mathbf{a}_{i}^{\rm T} \mathbf{x}^{\star}))$ for a fixed true parameter $\mathbf{x}^{\star} \in \mathbb{R}^{N_{\mathbf{x}}}$. The true parameter $\mathbf{x}^{\star}$ was computed once via a single draw from $\mathcal{N}(\mathbf{0}, \mathbf{I}_{N_{\mathbf{x}} \times N_{\mathbf{x}}})$. 

\begin{figure}[h]
  \begin{center}
    {\def\arraystretch{1.5}\tabcolsep=3pt
    \begin{tabular}{| c | c |}
    \hline
    \includegraphics[width=5.5cm]{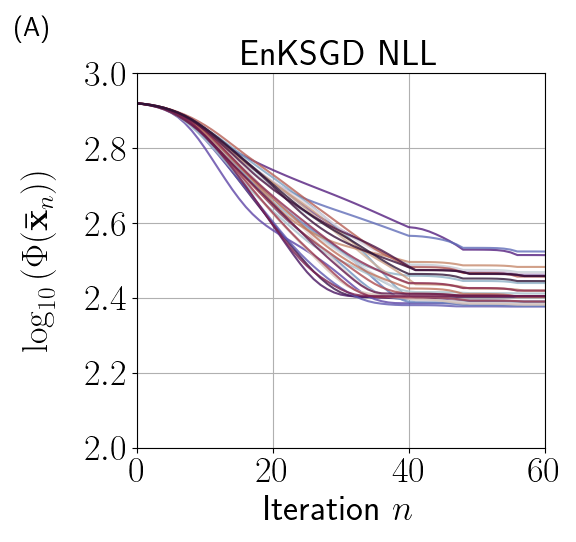} &
    \includegraphics[width=5.5cm]{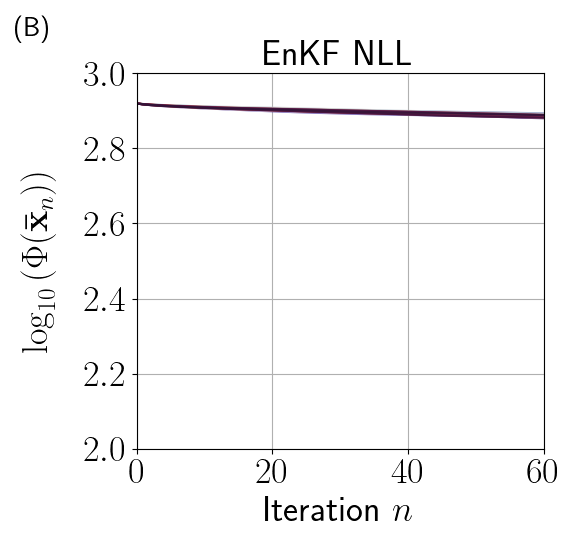} \\ 
    \hline
    \end{tabular}
    \caption{\textbf{Poisson Regression}: Results of $30$ independent runs of EnKSGD and EnKF when minimizing the NLL for Poisson regression with a simulated dataset. The average total number of $\mathcal{G}(\mathbf{x})$ evaluations for EnKSGD in (A) is $1585$, and for EnKF in (B) is $1561$. }
    \label{fig:poisson-regression-all}
    }
  \end{center}
\end{figure}

Given the simulated data $\mathbf{a}_{i}$ and $b_i$, Figure~\ref{fig:poisson-regression-all} compares EnKSGD vs. the EnKF-type approach when minimizing the Poisson regression NLL. The number of particles $K = 25$, the parameters $\beta = 10^{-6}$ and $\delta = 1$, and the initial mean $\mathbf{\bar{x}}_{0} = 2.5 \cdot \textbf{1}_{N_{\mathbf{x}}}$.

\subsection{Signal Reconstruction}
\label{subsec:signal-reconstruction}
Finally, we consider the nonlinear problem
\begin{equation} \label{eq:signal-reconstruction-inversion-objective}
\min_{\mathbf{x}} \bigg \{ \frac{1}{2} \norm{\mathcal{G}(\mathbf{x}) - \mathbf{y}_{\rm obs} }_2^2 + \frac{\alpha_{\mathbf{x}}}{2} \norm{\mathbf{W} \mathbf{x}}_2^2 + \frac{\alpha_{\mathbf{y}}}{2} \norm{\mathbf{F} \mathcal{G}(\mathbf{x})}_2^2 \bigg \} 
\end{equation}
where $N_{\mathbf{x}} = N_{\mathbf{y}}$, $\mathbf{W} \in \mathbb{R}^{N_{\mathbf{x}} \times N_{\mathbf{x}}}$ is zero except for a $1$ in position $1$ and position $N_{\mathbf{x}}$ along the main diagonal, and $\mathbf{F} \in \mathbb{R}^{(N_{\mathbf{x}}-1) \times N_{\mathbf{x}}}$ is a forward finite difference matrix which is zero except for $1$ along the main diagonal and $-1$ along the super-diagonal. The forward map $\big [ \mathcal{G}(\mathbf{x}) \big ]_{i} = 100 \tanh(\mathbf{x}_{i}/25), \forall i \in \{ 1, \dots, N_{\mathbf{y}} \}$, models a nonlinear amplifier. Minimizing (\ref{eq:signal-reconstruction-inversion-objective}) is one way to reconstruct/recover an input signal $\mathbf{x}(t)$ from noisy measurements $\mathbf{y}_{\rm obs}$ of the amplifier output. 

\begin{figure}[h]
  \begin{center}
    {\def\arraystretch{1.5}\tabcolsep=3pt
    \begin{tabular}{| c | c |}
    \hline
    \includegraphics[width=5.5cm]{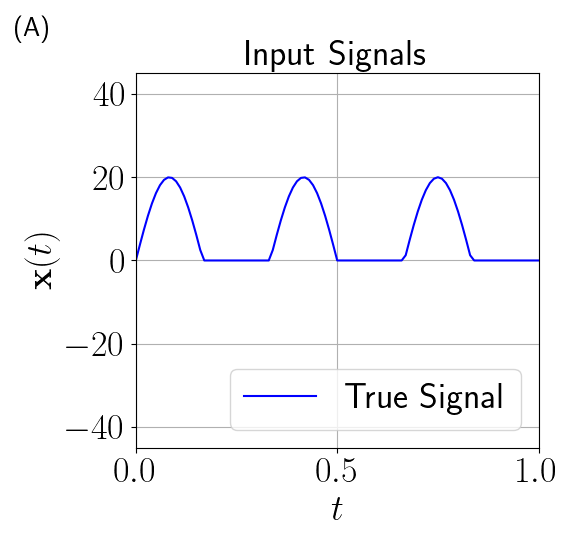} &
    \includegraphics[width=5.5cm]{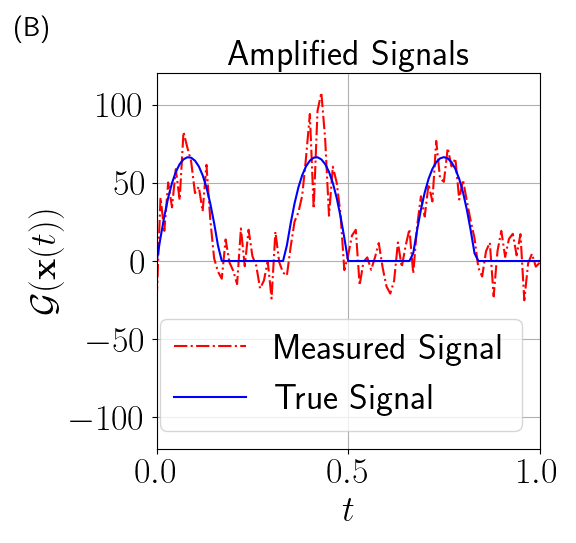} \\ 
    \hline
    \includegraphics[width=5.3cm, height=5.3cm]{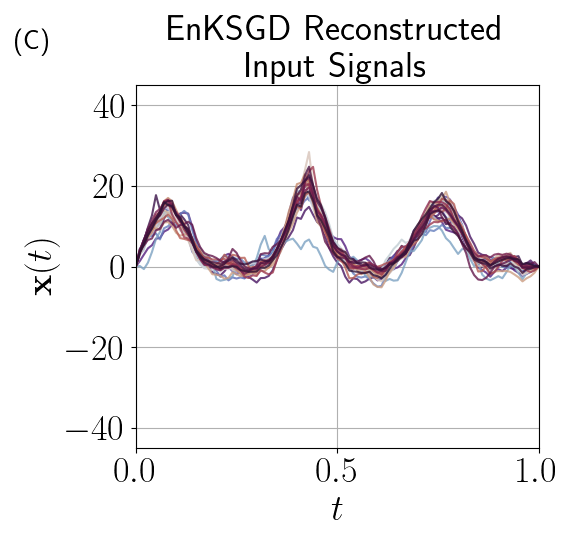} &
    \includegraphics[width=5.3cm, height=5.3cm]{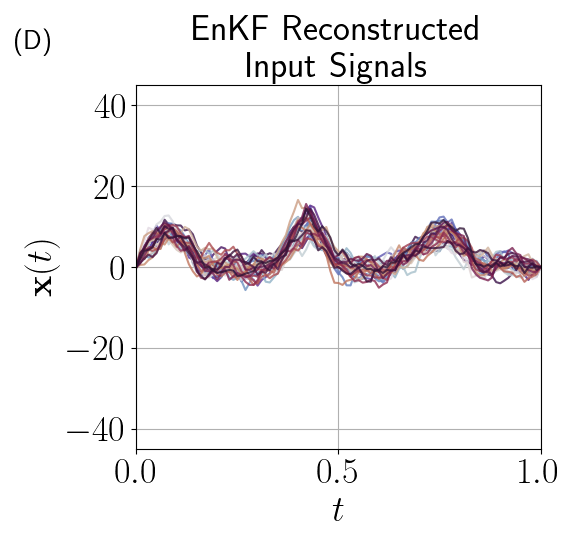} \\
    \hline
    \end{tabular}
    \caption{\textbf{Signal Reconstruction}: (C) and (D) each show $30$ independent reconstructions of (A) obtained by minimizing (\ref{eq:signal-reconstruction-inversion-objective}). After $60$ iterations of EnKSGD, the average total number of $\mathcal{G}(\mathbf{x})$ evaluations is $6195$ and the mean value of $\log_{10}(\Phi(\mathbf{\bar{x}}_{60}))$ is $4.28$. After $61$ iterations of EnKF, the average total number of $\mathcal{G}(\mathbf{x})$ evaluations is $6223$ and the mean value of $\log_{10}(\Phi(\mathbf{\bar{x}}_{61}))$ is $4.50$. }
    \label{fig:signal-reconstruction-all}
    }
  \end{center}
\end{figure}

Part (A) of Figure~\ref{fig:signal-reconstruction-all} shows the input signal $\mathbf{x}(t)$ we wish to recover, which is the half-wave rectified version of the sine wave $20 \sin(6 \pi t)$. Part (B) shows $\mathbf{x}(t)$ after amplification via $\mathcal{G}(\mathbf{x})$, and $N_{\mathbf{y}} = 101$ measurements $\mathbf{y}_{\rm obs}$ of the amplified signal corrupted by IID zero-mean Gaussian noise with a standard deviation of $15$. Parts (C) and (D) show reconstructions of $\mathbf{x}(t)$ obtained by minimizing (\ref{eq:signal-reconstruction-inversion-objective}) with $\alpha_{\mathbf{x}} = 10^{10}$ and $\alpha_{\mathbf{y}} = 5$. Note that increasing the value of $\alpha_{\mathbf{x}}$ forces the reconstructed signal towards zero at the first and last measurements, and increasing the value of $\alpha_{\mathbf{y}}$ increases the degree of smoothing. The results in Figure~\ref{fig:signal-reconstruction-all} were generated using $K = 101$ particles, the parameters $\beta = 10^{-6}$ and $\delta = 10^{-3}$, and the initial mean $\mathbf{\bar{x}}_{0} = \mathbf{0}$.

\subsection{Discussion}
\label{subsec:discussion-of-numerical-experiments}
Overall, the numerical results demonstrate that EnKSGD is robust to noise and can consistently outperform alternative approaches, sometimes by many orders of magnitude. In subsection~\ref{subsec:ill-conditioned-linear-least-squares}, EnKSGD outperforms both alternatives by over $10$ orders of magnitude, and comes by far the closest to the noise level of $\Phi(\mathbf{x})$. The CFD approach is not competitive with EnKSGD, and is most negatively impacted by noise. The CFD approach also exhibits the slowest inital decrease, although without noise it reaches a lower terminal value of $\Phi(\mathbf{x})$ than the EnKF-type approach. In subsection~\ref{subsec:nonlinear-least-squares}, EnKSGD outperforms the EnKF-type approach with respect to both the mean and median on $8/11 \approx 73 \%$ of the problems, and performs at least as well as the EnKF-type approach on all problems. Furthermore, EnKSGD beats the EnKF-type approach by $20$ orders of magnitude on the well-known Rosenbrock function. Observe that EnKSGD can outperform the EnKF-type approach when the number of particles $K > N_{\mathbf{x}}$ (e.g. subsections 5.1 and 5.2), when $K = N_{\mathbf{x}}$ (e.g. subsection 5.4), and when $K < N_{\mathbf{x}}$ (e.g. subsections 5.2 and 5.3).

\section{Conclusions}
\label{sec:conclusions}
In this paper, we introduced the EnKSGD class of algorithms. The EnKSGD class of algorithms improves upon previous EnKF approaches to optimization by: (i) accelerating convergence in the context of derivative-free optimization; (ii) enabling the use of a backtracking line search; and (iii) generalizing beyond the $L^{2}$ loss. Numerical experiments using linear least squares, nonlinear least squares, and maximum likelihood estimation problems demonstrate that EnKSGD is robust, fast, and applicable to a variety of important applied problems in scientific computing and machine learning. In future work, we believe it is worthwhile to: (i) study the performance of EnKSGD on empirical risk minimization problems with large datasets where mini-batching becomes necessary; (ii) investigate what convergence guarantees can be established for EnKSGD; and (iii) extend EnKSGD to handle constrained optimization problems.

\section*{Acknowledgments}
\label{sec:acknowledgments}
Brian Irwin thanks Eldad Haber for support and helpful discussions, and the members of SFB 1294 {\it Data Assimilation} for excellent working conditions at the University of Potsdam.

\bibliographystyle{siamplain}
\bibliography{references}

\begin{thebibliography}{10}

\bibitem{albers-EKI-constrained-2019}
{\sc D.~J. Albers, P.-A. Blancquart, M.~E. Levine, E.~Esmaeilzadeh~Seylabi, and
  A.~Stuart}, {\em Ensemble {K}alman methods with constraints}, Inverse
  problems, 35 (2019), p.~95007,
  \url{https://doi.org/10.1088/1361-6420/ab1c09}.

\bibitem{amezcua-ETKBF-2014}
{\sc J.~Amezcua, K.~Ide, E.~Kalnay, and S.~Reich}, {\em Ensemble transform
  {K}alman–{B}ucy filters}, Quarterly Journal of the Royal Meteorological
  Society, 140 (2014), pp.~995--1004, \url{https://doi.org/10.1002/qj.2186}.

\bibitem{asch-data-assimilation-book-2016}
{\sc M.~Asch, M.~Bocquet, and M.~Nodet}, {\em Data Assimilation}, Society for
  Industrial and Applied Mathematics, Philadelphia, PA, 2016,
  \url{https://doi.org/10.1137/1.9781611974546}.

\bibitem{bellantoni-Kalman-Schmidt-square-root-1967}
{\sc J.~F. Bellantoni and K.~W. Dodge}, {\em A square root formulation of the
  {K}alman-{S}chmidt filter}, AIAA journal, 5 (1967), pp.~1309--1314,
  \url{https://doi.org/10.2514/3.4189}.

\bibitem{bergemann-EnKF-localization-2010}
{\sc K.~Bergemann and S.~Reich}, {\em A localization technique for ensemble
  {K}alman filters}, Quarterly Journal of the Royal Meteorological Society, 136
  (2010), pp.~701--707, \url{https://doi.org/10.1002/qj.591}.

\bibitem{bergemann-mollified-EnKF-2010}
{\sc K.~Bergemann and S.~Reich}, {\em A mollified ensemble {K}alman filter},
  Quarterly Journal of the Royal Meteorological Society, 136 (2010),
  pp.~1636--1643, \url{https://doi.org/10.1002/qj.672}.

\bibitem{bezanson-julia-2017}
{\sc J.~Bezanson, A.~Edelman, S.~Karpinski, and V.~B. Shah}, {\em Julia: A
  fresh approach to numerical computing}, SIAM Review, 59 (2017), pp.~65--98,
  \url{https://doi.org/10.1137/141000671}.

\bibitem{bishop-ETKF-2001}
{\sc C.~H. Bishop, B.~J. Etherton, and S.~J. Majumdar}, {\em Adaptive sampling
  with the ensemble transform {K}alman filter. {P}art {I}: {T}heoretical
  aspects}, Monthly Weather Review, 129 (2001), pp.~420 -- 436,
  \url{https://doi.org/10.1175/1520-0493(2001)129<0420:ASWTET>2.0.CO;2}.

\bibitem{blomker-EKI-convergent-scheme-2018}
{\sc D.~Bl\"{o}mker, C.~Schillings, and P.~Wacker}, {\em A strongly convergent
  numerical scheme from ensemble {K}alman inversion}, SIAM Journal on Numerical
  Analysis, 56 (2018), pp.~2537--2562,
  \url{https://doi.org/10.1137/17M1132367}.

\bibitem{blomker-EKI-convergence-analysis-2019}
{\sc D.~Bl\"{o}mker, C.~Schillings, P.~Wacker, and S.~Weissmann}, {\em Well
  posedness and convergence analysis of the ensemble {K}alman inversion},
  Inverse Problems, 35 (2019), p.~085007,
  \url{https://doi.org/10.1088/1361-6420/ab149c}.

\bibitem{calvello-ensemble-Kalman-mean-field-perspective-ArXiv-2022}
{\sc E.~Calvello, S.~Reich, and A.~M. Stuart}, {\em Ensemble {K}alman methods:
  {A} mean field perspective},  (2022), \url{https://arxiv.org/abs/2209.11371}.

\bibitem{chada-EKI-parameterizations-2018}
{\sc N.~K. Chada, M.~A. Iglesias, L.~Roininen, and A.~M. Stuart}, {\em
  Parameterizations for ensemble {K}alman inversion}, Inverse Problems, 34
  (2018), p.~055009, \url{https://doi.org/10.1088/1361-6420/aab6d9}.

\bibitem{chada-EKI-tikhonov-2020}
{\sc N.~K. Chada, A.~M. Stuart, and X.~T. Tong}, {\em Tikhonov regularization
  within ensemble {K}alman inversion}, SIAM Journal on Numerical Analysis, 58
  (2020), pp.~1263--1294, \url{https://doi.org/10.1137/19M1242331}.

\bibitem{ding-EKI-analysis-2021}
{\sc Z.~Ding and Q.~Li}, {\em Ensemble {K}alman inversion: {M}ean-field limit
  and convergence analysis}, Statistics and computing, 31 (2021),
  \url{https://doi.org/10.1007/s11222-020-09976-0}.

\bibitem{ding-EKS-analysis-2021}
{\sc Z.~Ding and Q.~Li}, {\em Ensemble {K}alman sampler: {M}ean-field limit and
  convergence analysis}, SIAM Journal on Mathematical Analysis, 53 (2021),
  pp.~1546--1578, \url{https://doi.org/10.1137/20M1339507}.

\bibitem{ernst-EnKF-analysis-2015}
{\sc O.~G. Ernst, B.~Sprungk, and H.-J. Starkloff}, {\em Analysis of the
  ensemble and polynomial chaos {K}alman filters in {B}ayesian inverse
  problems}, SIAM/ASA Journal on Uncertainty Quantification, 3 (2015),
  pp.~823--851, \url{https://doi.org/10.1137/140981319}.

\bibitem{evensen-EnKF-original-paper-1994}
{\sc G.~Evensen}, {\em Sequential data assimilation with a nonlinear
  quasi-geostrophic model using {M}onte {C}arlo methods to forecast error
  statistics}, Journal of Geophysical Research: Oceans, 99 (1994),
  pp.~10143--10162, \url{https://doi.org/10.1029/94JC00572}.

\bibitem{evensen-EnKF-formulation-and-implementation-2003}
{\sc G.~Evensen}, {\em The ensemble {K}alman filter: {T}heoretical formulation
  and practical implementation}, Ocean dynamics, 53 (2003), pp.~343--367,
  \url{https://doi.org/10.1007/s10236-003-0036-9}.

\bibitem{evensen-EnKF-book-2009}
{\sc G.~Evensen}, {\em Data assimilation: {T}he ensemble {K}alman filter},
  Springer, 2nd~ed., 2009.

\bibitem{garbuno-inigo-EKS-2020}
{\sc A.~Garbuno-Inigo, F.~Hoffmann, W.~Li, and A.~M. Stuart}, {\em Interacting
  {L}angevin diffusions: {G}radient structure and ensemble {K}alman sampler},
  SIAM Journal on Applied Dynamical Systems, 19 (2020), pp.~412--441,
  \url{https://doi.org/10.1137/19M1251655}.

\bibitem{haber-never-look-back-EnKF-ArXiv-2018}
{\sc E.~Haber, F.~Lucka, and L.~Ruthotto}, {\em Never look back - {A} modified
  {E}n{KF} method and its application to the training of neural networks
  without back propagation},  (2018), \url{https://arxiv.org/abs/1805.08034}.

\bibitem{huang-derivative-free-Bayesian-inversion-2022}
{\sc D.~Z. Huang, J.~Huang, S.~Reich, and A.~M. Stuart}, {\em Efficient
  derivative-free {B}ayesian inference for large-scale inverse problems},
  Inverse problems, 38 (2022), \url{https://doi.org/10.1088/1361-6420/ac99fa}.

\bibitem{huang-iterated-Kalman-methodology-2022}
{\sc D.~Z. Huang, T.~Schneider, and A.~M. Stuart}, {\em Iterated {{Kalman}}
  methodology for inverse problems}, Journal of Computational Physics, 463
  (2022), p.~111262, \url{https://doi.org/10.1016/j.jcp.2022.111262}.

\bibitem{iglesias-EKI-pde-constrained-2016}
{\sc M.~A. Iglesias}, {\em A regularizing iterative ensemble {K}alman method
  for {PDE}-constrained inverse problems}, Inverse problems, 32 (2016),
  pp.~25002--25046, \url{https://doi.org/10.1088/0266-5611/32/2/025002}.

\bibitem{iglesias-EKI-2013}
{\sc M.~A. Iglesias, K.~J.~H. Law, and A.~M. Stuart}, {\em Ensemble {K}alman
  methods for inverse problems}, Inverse problems, 29 (2013), pp.~45001--20,
  \url{https://doi.org/10.1088/0266-5611/29/4/045001}.

\bibitem{julier-UKF-1997}
{\sc S.~J. Julier and J.~K. Uhlmann}, {\em New extension of the {K}alman filter
  to nonlinear systems}, in Signal Processing, Sensor Fusion, and Target
  Recognition VI, vol.~3068, Bellingham WA, 1997, SPIE, pp.~182--193,
  \url{https://doi.org/10.1117/12.280797}.

\bibitem{kalman-new-linear-filtering-approach-1960}
{\sc R.~E. Kalman}, {\em A new approach to linear filtering and prediction
  problems}, Journal of basic engineering, 82 (1960), pp.~35--45,
  \url{https://doi.org/10.1115/1.3662552}.

\bibitem{kelly-EnKF-well-posedness-2014}
{\sc D.~T.~B. Kelly, K.~J.~H. Law, and A.~M. Stuart}, {\em Well-posedness and
  accuracy of the ensemble {K}alman filter in discrete and continuous time},
  Nonlinearity, 27 (2014), pp.~2579--2603,
  \url{https://doi.org/10.1088/0951-7715/27/10/2579}.

\bibitem{kovachki-EKI-machine-learning-2019}
{\sc N.~B. Kovachki and A.~M. Stuart}, {\em Ensemble {K}alman inversion: {A}
  derivative-free technique for machine learning tasks}, Inverse problems, 35
  (2019), p.~95005, \url{https://doi.org/10.1088/1361-6420/ab1c3a}.

\bibitem{kwiatkowski-square-root-EnKF-convergence-2015}
{\sc E.~Kwiatkowski and J.~Mandel}, {\em Convergence of the square root
  ensemble {K}alman filter in the large ensemble limit}, SIAM/ASA Journal on
  Uncertainty Quantification, 3 (2015), pp.~1--17,
  \url{https://doi.org/10.1137/140965363}.

\bibitem{larson-dfo-review-2019}
{\sc J.~Larson, M.~Menickelly, and S.~M. Wild}, {\em Derivative-free
  optimization methods}, Acta Numerica, 28 (2019), pp.~287--404,
  \url{https://doi.org/10.1017/S0962492919000060}.

\bibitem{law-data-assimilation-book-2015}
{\sc K.~Law, A.~Stuart, and K.~Zygalakis}, {\em Data Assimilation: {A}
  Mathematical Introduction}, Springer International Publishing, 2015,
  \url{https://doi.org/10.1007/978-3-319-20325-6}.

\bibitem{law-deterministic-EnKF-2016}
{\sc K.~J.~H. Law, H.~Tembine, and R.~Tempone}, {\em Deterministic mean-field
  ensemble {K}alman filtering}, SIAM Journal on Scientific Computing, 38
  (2016), pp.~A1251--A1279, \url{https://doi.org/10.1137/140984415}.

\bibitem{legland-EnKF-asymptotics-2009}
{\sc F.~Le~Gland, V.~Monbet, and V.-D. Tran}, {\em Large sample asymptotics for
  the ensemble {K}alman filter}, Tech. Report RR-7014, INRIA, 2009,
  \url{https://hal.inria.fr/inria-00409060}.

\bibitem{li-EnKF-numerical-properties-2008}
{\sc J.~Li and D.~Xiu}, {\em On numerical properties of the ensemble {K}alman
  filter for data assimilation}, Computer methods in applied mechanics and
  engineering, 197 (2008), pp.~3574--3583,
  \url{https://doi.org/10.1016/j.cma.2008.03.022}.

\bibitem{nocedal-wright-numerical-optimization-2006}
{\sc J.~Nocedal and S.~J. Wright}, {\em Numerical Optimization},
  Springer-Verlag, New York, 2006,
  \url{https://doi.org/10.1007/978-0-387-40065-5}.

\bibitem{nusken-EKS-note-2019}
{\sc N.~N\"{u}sken and S.~Reich}, {\em Note on interacting {L}angevin
  diffusions: {G}radient structure and ensemble {K}alman sampler by
  {G}arbuno-{I}nigo, {H}offmann, {L}i and {S}tuart},  (2019),
  \url{https://arxiv.org/abs/1908.10890}.

\bibitem{opper-bayesian-online-learning-1999}
{\sc M.~Opper}, {\em A {B}ayesian Approach to On-line Learning}, Publications
  of the Newton Institute, Cambridge University Press, 1999, p.~363–378,
  \url{https://doi.org/10.1017/CBO9780511569920.017}.

\bibitem{orban-siqueira-nlsproblems-2021}
{\sc D.~Orban, A.~S. Siqueira, and {contributors}}, {\em {NLSProblems.jl}:
  Nonlinear least-squares problems for {NLPModels}}.
\newblock \url{https://github.com/JuliaSmoothOptimizers/NLSProblems.jl}, March
  2021, \url{https://doi.org/10.5281/zenodo.4605405}.

\bibitem{pidstrigach-ensemble-transform-logistic-2022}
{\sc J.~Pidstrigach and S.~Reich}, {\em Affine-invariant ensemble transform
  methods for logistic regression}, Foundations of computational mathematics,
  (2022), \url{https://doi.org/10.1007/s10208-022-09550-2}.

\bibitem{pinnau-consensus-optimization-2017}
{\sc R.~Pinnau, C.~Totzeck, O.~Tse, and S.~Martin}, {\em A consensus-based
  model for global optimization and its mean-field limit}, Mathematical Models
  and Methods in Applied Sciences, 27 (2017), pp.~183--204,
  \url{https://doi.org/10.1142/S0218202517400061}.

\bibitem{reich-dynamical-2011}
{\sc S.~Reich}, {\em A dynamical systems framework for intermittent data
  assimilation}, BIT Numerical Mathematics, 51 (2011), pp.~235--249,
  \url{https://doi.org/10.1007/s10543-010-0302-4}.

\bibitem{reich-cotter-book-2015}
{\sc S.~Reich and C.~Cotter}, {\em Probabilistic forecasting and {B}ayesian
  data assimilation}, Cambridge University Press, Cambridge, 2015,
  \url{https://doi.org/10.1017/CBO9781107706804}.

\bibitem{schillings-EnKF-analysis-2017}
{\sc C.~Schillings and A.~M. Stuart}, {\em Analysis of the ensemble {K}alman
  filter for inverse problems}, SIAM Journal on Numerical Analysis, 55 (2017),
  pp.~1264--1290, \url{https://doi.org/10.1137/16M105959X}.

\bibitem{schillings-EKI-convergence-linear-noisy-2018}
{\sc C.~Schillings and A.~M. Stuart}, {\em Convergence analysis of ensemble
  {K}alman inversion: {T}he linear, noisy case}, Applicable Analysis, 97
  (2018), pp.~107--123, \url{https://doi.org/10.1080/00036811.2017.1386784}.

\bibitem{tong-nonlinear-stability-EnKF-2016}
{\sc X.~T. Tong, A.~J. Majda, and D.~Kelly}, {\em Nonlinear stability and
  ergodicity of ensemble based {K}alman filters}, Nonlinearity, 29 (2016),
  pp.~657--691, \url{https://doi.org/10.1088/0951-7715/29/2/657}.

\end{thebibliography}
\end{document}